\documentclass[preprint,12pt]{article}
\usepackage{amssymb,amsmath,hyperref,enumerate}
\usepackage[english]{babel}

\usepackage[all]{xy}

\usepackage{amssymb}

\newtheorem{thm}{Theorem}[section]
\newtheorem{Theorem}{Main Result}
\newtheorem{rem}[thm]{Remark}
\newtheorem{prop}[thm]{Proposition}
\newtheorem{lemma}[thm]{Lemma}
\newtheorem{cor}[thm]{Corollary}
\newtheorem{Cor}{Main Corollary}

\def\<{\langle}
\def\>{\rangle}
\newcommand{\proof}{\emph{Proof.~}}
\newcommand{\cC}{\mathcal{C}}
\newcommand{\Aut}{\mathrm{Aut}}

\newcommand{\cP}{\mathcal{P}}

\def\qed{{\hfill\hphantom{.}\nobreak\hfill$\Box$}}
\newcommand{\proj}{\mathrm{proj}}
\newcommand{\cL}{\mathcal{L}}

\newcommand{\id}{\mathrm{id}}

\newcommand{\R}{\mathbb{R}}
\newcommand{\N}{\mathbb{N}}

\newcommand{\inc}{\mbox{\tt I}}

\begin{document}




\title{On epimorphisms of spherical Moufang buildings}


\author{Koen Struyve\thanks{The author is supported by the Fund for Scientific Research --
Flanders (FWO - Vlaanderen)}}
%

\maketitle
\begin{abstract}
In this paper we classify the the epimorphisms of irreducible spherical Moufang buildings (of rank $\geq 2$) defined over a field. As an application we characterize indecomposable epimorphisms of these buildings as those epimorphisms arising from $\R$-buildings. 
\end{abstract}








\section{Introduction}

The theory of buildings was introduced by Jacques Tits in the late 60's in order to better understand certain classes of (algebraic) groups. This theory certainly attained this goal and much more. The two most studied subclasses are the spherical and affine buildings.

The spherical buildings have been classified by Jacques Tits in 1974 (\cite{Tit:74}) provided that the rank is at least three. Using the so-called `spherical building at infinity' of an affine building, Tits also classified in the affine buildings of rank at least 4 (\cite{Tit:86}). This classification also includes non-discrete generalizations of affine buildings, the $\R$-buildings.

Whereas this classification uses spherical buildings to say something about $\R$-buildings, in the current paper we will use $\R$-buildings to answer a problem concerning spherical buildings. The question is to classify or characterize epimorphisms of Moufang spherical buildings. We will show that epimorphisms of Moufang buildings  defined over a field (made precise in Section~\ref{section:deffield})  correspond to valuations of the underlying field satisfying certain compatibility conditions. 
Epimorphisms arising from Moufang $\R$-buildings turn out to be the `primitive' epimorphisms for this class, i.e. if one cannot decompose the epimorphism into two proper epimorphisms, then the epimorphism arises directly from an $\R$-building.

For a precise version of the main results and corollaries we refer to Section~\ref{section:main}.

This extends known results for projective spaces (see Section~\ref{section:known}). The remaining open class, consisting of certain polar spaces of pseudo-quadratic form type defined over a nonabelian skew field is handled by the author and Petra N. Schwer in a forthcoming paper (\cite{Sch-Str:*}) using different, case-specific methods.

\textbf{Acknowledgement.} The author would like to thank Pierre-Emmanuel Caprace for suggesting the problem.

\section{Preliminaries}

\subsection{Buildings}\label{section:buildings}
Let $(W,S)$ be a Coxeter system, then a \emph{weak building of type $(W,S)$} is a pair $(\mathcal{C},\delta)$ consisting of a nonempty set $\mathcal{C}$ (called \emph{chambers}) and a map $\delta:\mathcal{C} \times \mathcal{C} \rightarrow W$ (called the \emph{Weyl distance}), such that for every two chambers $C$ and $D$ the following holds.
\begin{itemize}
\item[{(WD1)}]
 $\delta(C,D) = 1$ if and only if $C=D$.
\item[{(WD2)}]
If $\delta(C,D) = w$ and $C' \in \mathcal{C}$ satisfies $\delta(C',C)=s \in S$, then $\delta(C',D) \in \{sw,w\}$. If moreover $l(sw)=l(w)+1$ (where $l$ is the word metric on $W$ w.r.t. $S$), then $\delta(C',D) =sw$.
\item[{(WD3)}]
If $\delta(C,D) = w$, then for any $s \in S$ there exists a chamber $C' \in \mathcal{C}$ such that $\delta'(C',C) =s$ and $\delta(C',D)=sw$. 
\end{itemize}

This weak building is said to be \emph{spherical} if the Coxeter group $W$ is finite. The \emph{rank} of a weak building is defined to be $\vert S \vert$. Two chambers are \emph{$s$-equivalent} (with $s\in S$) if the Weyl distance between them is either $s$ or the identity element $1$ of $W$. Consider a subset $S' \subset S$. The connected components of $\cC$ using only equivalences in $S'$ are called the \emph{$S'$-residues}, which are again weak buildings. The rank one residues (so $S'=\{s\}$) are also called \emph{($s$-)panels}. If each panel of the weak building has cardinality at least 3, then we say that $(\mathcal{C},\delta)$ is a \emph{building}. 


A building is \emph{irreducible} if it cannot be decomposed as a direct product of two (non-trivial) buildings. 

A \emph{morphism} $\phi$ between two (weak) buildings $(\mathcal{C},\delta)$ and $(\mathcal{C}',\delta')$ of type $(W,S)$ is a map from $\mathcal{C}$ to $\mathcal{C}'$ preserving $s$-equivalency for each $s\in S$. If in addition this map is respectively injective or surjective, then it is respectively called an \emph{endomorphism} or an \emph{epimorphism}. If it is both injective and surjective then it is an \emph{isomorphism}. We say that an automorphism $g$ of a building $(\mathcal{C},\delta)$ \emph{descends} under an epimorphism $\phi$ from $(\mathcal{C},\delta)$ to $(\mathcal{C}',\delta')$ if there exists an automorphism $g'$ of $(\mathcal{C}',\delta')$ such that $\phi \circ g = g' \circ \phi$. One easily verifies this is equivalent to the condition $\forall C,D \in \mathcal{C}: C^\phi = D^\phi \Leftrightarrow C^{g \phi} = D^{g \phi}$. Note that the automorphisms that descend form a subgroup of the full automorphism group.

Two chambers of a building of spherical building are \emph{opposite} if the Weyl distance between them is maximal with respect to the word metric on $W$. An $s$-panel and $s'$-panel are \emph{opposite} if $s$ and $s'$ are mapped to each other by the opposition involution of the Coxeter group (see~\cite[p. 61]{Abr-Bro:08}) and contain opposite chambers. Opposite panels have the property that for each chamber in one of these panels there is a unique non-opposite chamber in the other one. 

For more information on (spherical) buildings, we refer to~\cite{Abr-Bro:08} and~\cite{Wei:03}.

\begin{rem} \rm
If we speak about an epimorphism of a building, we assume that its image is not a weak building. Non-type preserving epimorphisms will not be considered in this paper.
\end{rem}

\begin{rem} \rm
Our main results only deal with irreducible buildings. As reducible buildings are direct products of irreducible buildings, the study of the epimorphisms of these can be reduced to their components. 
\end{rem}


\subsection{Generalized polygons}
For the spherical buildings of rank 2, the generalized polygons, we will take an incidence geometric point of view using the panels as basic objects. We define them as follows.

Let $\Gamma:=(\cP,\cL,\inc)$ be a rank 2 geometry consisting of a \emph{point set} $\cP$, a \emph{line set} $\cL$ (with $\cP \cap \cL=\emptyset$), and \emph{incidence relation} $\inc$ between $\cP$ and $\cL$. An \emph{element} of $\Gamma$ is a point or line of it. An ordered sequence $(x_0, x_1, \dots, x_k)$ of elements of $\Gamma$ is called a \emph{path} of {\em length} $k$ if each two subsequent elements in it are incident. We say it \emph{stammers} if there is an $i$ such that $x_i$ and $x_{i+2}$ are identical. 

The rank 2 geometry $\Gamma=(\cP,\cL,\inc)$ is a generalized $n$-gon ($n \in \N$, $n\geq2$) if it satisfies the following axioms.
\begin{itemize}
 \item[(GP1)] Every element is incident with at least three other elements.
 \item[(GP2)] For every pair of elements $x,y\in \cP\cup \cL$,
 there exists a non-stammering path $(x_0=x,x_1,\ldots,x_{k-1},x_k=y)$ of length at most $n$
 \item[(GP3)] The sequence in (GP2) is unique if its length is strictly smaller than $n$.
\end{itemize}

Note that this definition is self-dual in the notions point and line. The chambers here are incident point-line pairs. Panels are the sets of chambers containing a certain element. The corresponding building is irreducible if and only if $n\geq 3$. We define an \emph{apartment} to be a substructure that is an ordinary $n$-gon. A \emph{root} is a non-stammering path of length $n$.


The \emph{distance} between two elements is the length of a shortest path between them. Two elements at maximal distance $n$ are said to be \emph{opposite}. If $x$ and $y$ are not opposite or equal, then the \emph{projection of $y$ on $x$} (denoted by $\proj_x y$) is the unique element incident with $x$ closest to $y$. 

Morphisms from this point of view are maps between generalized $n$-gons, mapping points to points, lines to lines, such that incident elements are mapped to incident elements. Endomorphisms and epimorphisms are then defined as usual.

If the image of a non-stammering path under a epimorphism of the generalized polygon becomes stammering, we say that the path \emph{collapses} under the epimorphism. We will use the same notion for apartments, by considering a non-stammering path of length $2n$ defining the apartment.

\subsubsection{The Moufang property}

Let $\alpha:= (x_0, x_1, \dots, x_n)$ be a root of a generalized $n$-gon $\Gamma$ with $n\geq 3$. A \emph{root elation
 of $\alpha$} is an automorphism of $\Gamma$ fixing each element incident with an element of the subpath $(x_1, \dots, x_{n-1})$. The \emph{root group of $\alpha$} is the group consisting of all root elations of $\alpha$. We say that $\alpha$ is  \emph{Moufang} if this group acts transitively on the elements incident with $x_0$ different from $x_1$. One shows that if this is the case then the group acts sharply transitive on this set (see for instance~\cite[Def. 5.2.1]{Mal:98}). The generalized polygon $\Gamma$ is \emph{Moufang} if all its roots are.
 
\begin{rem} \rm
It is possible to generalize this definition to higher rank (spherical) buildings. An irreducible spherical building of rank at least 3 is automatically Moufang by a result of Tits (\cite{Tit:74}). 
\end{rem}

\subsection{Classifications of spherical buildings and the underlying field}
The book~\cite{Tit:74} of Tits includes a classification of the irreducible spherical buildings of rank at least 3. Moufang generalized polygons have been classified by Tits and Weiss in~\cite{Tit-Wei:02}.

The aim of this section is to briefly discuss this classification and clarify what we mean by `defined over a field' and `underlying field' in the statement of the main results and corollaries (Section~\ref{section:main}). These notions are not unambiguous and will be different than the point of view of~\cite[Rem. 30.29]{Wei:09}.

\subsubsection{Moufang generalized polygons}\label{section:Moufang}
We start with the Moufang generalized polygons (in which we follow~\cite{Tit-Wei:02}). Let $\Sigma$ be an apartment of a generalized Moufang $n$-gon $\Gamma$, and label the elements of it by $x_i$, with $i \in \mathbb{Z}$ such that $x_i \inc x_{i+1}$ and $x_i = x_{i+2n}$. This apartment will be called the \emph{hat-rack}. Let $U_i$ be the root group of the root $(x_i, x_{i+1}, \dots , x_{i+n})$. All of the $U_i$ forms the \emph{root group data} of $\Gamma$ associated to $\Sigma$. We will often use subscripts to indicate to which root group an automorphism belongs. 

Define $U_{[i,j]}$ to be the group generated by $U_i, U_{i+1}, \dots ,U_j$ (if $j<i$, then we let $U_{[i,j]}$ denote the group consisting only of the trivial automorphism).  

The following lemmas express the commutator relations between $U_i$ and $U_j$ when the corresponding roots are not opposite (i.e. $i \not\equiv j \mod 2n$).

\begin{lemma}[\cite{Tit-Wei:02}, Prop. 5.5] \label{lemma:nil1}
If $i+1 \leq j \leq i+n-1$, then $[U_i, U_j] \leq U_{[i+1,j-1]}$. \qed
\end{lemma}

\begin{lemma}[\cite{Tit-Wei:02}, Prop. 5.6] \label{lemma:nil2}
If $i \leq j \leq i+n -1$, then the product $U_i U_{i+1} \dots U_j$ is the group $U_{[i,j]}$, and every element of this group has a unique decomposition as $u_i u_{i+1} \dots u_j$ with $u_k \in U_k$. \qed
\end{lemma}

The last lemma implies that by giving descriptions of the root groups $U_1$ up to $U_n$ and the commutator relations between them, that one can completely describe the group $U_{[1,n]}$. One calls the $(n+1)$-tuple $(U_{[1,n]}, U_1, \dots, U_n)$ a \emph{root group sequence} of which $U_1, \dots, U_n$ are the \emph{terms}. This root group sequence suffices to describe the Moufang generalized polygon up to isomorphism (see~\cite[Prop. 12.2-3]{Wei:03}).

This reduces the classification to determining the possible root groups $U_1, \dots, U_n$ and their commutator relations. Let us briefly list the possibilities (for a detailed description see~\cite[16.1-9]{Tit-Wei:02}).

\begin{itemize}
\item{The triangles $\mathcal{T}(K)$}. 
\item{The quadrangles $\mathcal{Q}_\mathcal{I}(K,K_0, \sigma)$ of involutory type.}
\item{The quadrangles $\mathcal{Q}_\mathcal{Q}(K,L_0, q)$ of quadratic form type.}
\item{The quadrangles $\mathcal{Q}_\mathcal{D}(K,K_0, L_0)$ of indifferent type.}
\item{The quadrangles $\mathcal{Q}_\mathcal{P}(K,K_0, \sigma, L_0,q)$ of pseudo-quadratic form type.}
\item{The quadrangles $\mathcal{Q}_\mathcal{E}(K,L_0, q)$ of type $\mathsf{E}_i$ ($i =6,7,8$).}
\item{The quadrangles $\mathcal{Q}_\mathcal{F}(K,L_0, q)$ of type $\mathsf{F}_4$.}
\item{The hexagons $\mathcal{H}(J,K,\#)$.}
\item{The octagons $\mathcal{O}(K,\sigma)$.}
\end{itemize}

For the remainder of this paper we will consider the quadrangles of involutory type to be the subclass of those of pseudo-quadratic form type for which the pseudo-quadratic form is 0-dimensional. However we will need the class of \emph{quadrangles of quadratic and honorary involutory type}. These are quadrangles of quadratic form type where the vector space $L_0$ over $K$ with quadratic form $q$ can be interpreted as a composition algebra over $K$ with norm $q$.  These can also interpreted as involutory quadrangles except when this composition algebra is an octonion algebra (in which case one calls them honorary). More information on this classification can be found in~\cite[(38.9)]{Tit-Wei:02} and~\cite[30.14-30.33]{Wei:09}.

We set the \emph{underlying skew field or octonion algebra} for all of these cases to be $K$. We consider quadrangles of quadratic and honorary involutory type to be of quadratic form type, so they are defined over the underlying field $K$, not over the composition algebra.

\subsubsection{Higher rank buildings}\label{section:high}

In order to describe the higher rank case as well, one considers the following reduction. Let $\Delta := (\cC, \delta)$ be an irreducible spherical Moufang building of type $(W,S)$. If $C$ is a  chamber in  $\cC$ then one can describe the building $\Delta$ by means of a root group labeling, which is a triple $(u,\Theta,\theta)$, where 
\begin{itemize}
\item
$u$ maps involutions $s$ to groups $u(s)$, which correspond to the panels containing $C$,
\item
$\Theta_{s,t}$ is a root group sequence of a generalized polygon for every ordered pair $(s,t)$ of non-commuting involutions in $S$, these correspond to the rank two residues containing $C$ which are not digons.
\item
$\theta_{s,t}$ is an isomorphism from $u(s)$ to the first term of $\Theta_{s,t}$ for each ordered pair $(s,t)$ of non-commuting involutions in $S$. 
\end{itemize}
This data determines the building up to isomorphism, see~\cite[Prop. 12.2-3]{Wei:03}. For the rank two case the root group labeling is essentially the root group sequence.


Let us list, without much detail, the possibilities with rank at least three (after~\cite[12.13-12.19]{Wei:03}), with as modifications considering involutory type as a subclass of pseudo-quadratic form type, except for those of quadratic or honorary type which we put in a completely separate class). We also list each time the different isomorphism classes of rank 2 residues which occur (apart from digons).
\begin{itemize}
\item $\mathsf{A}_l(K)$ ($l \geq 3$): $\mathcal{T}(K)$.
\item $\mathsf{B}_l(K, L_0, q)$ ($l \geq 3$): $\mathcal{T}(K)$, $\mathcal{Q}_\mathcal{Q}(K,L_0, q)$.
\item $\mathsf{C}_l(K, K_0, \sigma)$ of quadratic ($l \geq 3$) or honorary type ($l =3$): $\mathcal{T}(K)$, $\mathcal{Q}_\mathcal{Q}(K_0, K, N)$ (where $N$ is the norm induced on the composition algebra $K$ over $K_0$).
\item $\mathsf{BC}_l(K, K_0, \sigma, L_0, q)$ (except from those isomorphic to a $\mathsf{C}_l(K, K_0, \sigma)$ of quadratic type)  ($l \geq 3$): $\mathcal{T}(K)$, $\mathcal{Q}_\mathcal{P}(K,K_0, \sigma, L_0,q)$.
\item $\mathsf{D}_l(K)$ ($l \geq 4$): $\mathcal{T}(K)$.
\item $\mathsf{E}_l(K)$ ($l=6,7,8$): $\mathcal{T}(K)$.
\item $\mathsf{F}_4(K,F,\sigma)$: $\mathcal{T}(K)$, $\mathcal{T}(F)$, $\mathcal{Q}_\mathcal{Q}(F,K, N)$ (where $N$ is the norm induced on the composition algebra $K$ over $F$).
\end{itemize}

The first four classes can be considered as continuations of rank 2 cases (see the last rank 2 residue listed each time) and are hence defined for $l=2$ as well.

We define the \emph{underlying skew field or octonion algebra} for the first, second, fourth, fifth and sixth case to be $K$.  For the third case we set it to be $K_0$ and for the seventh case $F$. 


\subsubsection{Spherical Moufang buildings defined over a field}\label{section:deffield}
We now list the irreducible spherical Moufang buildings of rank at least two that we consider to be defined over a field. These are exactly the cases for which our main results hold.

\begin{itemize}
\item $\mathsf{A}_l(K)$ ($l \geq 2$) where $K$ is a field,
\item $ \mathsf{B}_l(K, L_0, q)$ ($l \geq 2$),
\item $\mathsf{C}_l(K, K_0, \sigma)$ of quadratic ($l \geq 2$) or honorary type ($l =2$ or $3$) ,

\item $\mathsf{BC}_l(K, K_0, \sigma, L_0, q)$: $\mathcal{T}(K)$, $\mathcal{Q}_\mathcal{P}(K,K_0, \sigma, L_0,q)$  where $K$ is a field ($l \geq 2$),
\item $\mathsf{D}_l(K)$ ($l \geq 2$),
\item $\mathsf{E}_l(K)$ ($l=6,7,8$),
\item $\mathsf{F}_4(K,F,\sigma)$,
\item The quadrangles $\mathcal{Q}_\mathcal{D}(K,K_0, L_0)$, $\mathcal{Q}_\mathcal{E}(K,L_0, q)$, and $\mathcal{Q}_\mathcal{F}(K,L_0, q)$,
\item The hexagons $\mathcal{H}(J,K,\#)$,
\item The octagons $\mathcal{O}(K,\sigma)$.
\end{itemize}

Note that a spherical Moufang building defined over a field might have rank 2 residues not defined over a field. This occurs for the case $\mathsf{C}_l(K, K_0, \sigma)$ ($l \geq 3$) of quadratic type with $K$ a quaternion algebra or of honorary type ($l=3$), and for $\mathsf{F}_4(K,F,\sigma)$, when $K$ is a quaternion or octonion algebra.

The only spherical Moufang buildings of rank at least two not defined over field (and hence not covered by the results of this paper) are the projective spaces $\mathsf{A}_l(K)$ where $K$ is a nonabelian skew field or octonion algebra, and those polar spaces $\mathsf{BC}_l(K, K_0, \sigma, L_0, q)$, where $K$ is a nonabelian skew field, not isomorphic to a polar space of type $\mathsf{C}_l(K, K_0, \sigma)$ of quadratic type.

\subsection{Root group labelings with epimorphism data}\label{section:rgl}

We now define a way to denote the algebraic information we will derive from a given epimorphism of an irreducible spherical Moufang building of rank at least two defined over a field (so one of those listed in Section~\ref{section:deffield}). A 7-tuple $(u,\Theta,\theta,K, \Lambda, \nu, k)$ will be a \emph{root group labeling with epimorphism data} if:
\begin{itemize}
\item
$(u, \Theta, \theta)$ is a root group labeling of a building $\Delta$ listed in Section~\ref{section:deffield}.  We moreover demand that this root group labeling is of the standard form given in~\cite[12.13-12.19]{Wei:03} for those of rank at least three, or arises from the root group sequences of standard form given in \cite[16.1-9]{Tit-Wei:02}. In particular we want that the  groups $v(i)$ equal the parameter groups, which are the domains of the parametrizing maps $t \to x_i(t)$ of the root groups, cfr.~\cite[12.12]{Wei:03}. 
\item
The $K$ denotes the underlying field of $\Delta$, and $\nu$ is a (possibly trivial) surjective valuation of this field to an ordered abelian group $\Lambda$. 
\item
For the rank two case, so a Moufang $n$-gon, one can define maps from the parameter groups of the root groups $U_1$ and $U_n$ to the underlying field $K$ in the following way (using the aforementioned standard form).
\begin{center}
\[ \begin{array}{lccc} 
& & \underline{i =1} &  \underline{i =n }\\
\\
\underline{n=3} & \mathcal{T}: & \id  & \id\\ 
\underline{n=4} & \mathcal{Q_Q}: & \id & q \\ 
& \mathcal{Q_D}: & a \mapsto a & a \mapsto a^2 \\ 
 & \mathcal{Q_P}: & (a,t) \mapsto t  & \id \\ 
 & \mathcal{Q_E}: & (a,t) \mapsto q(\pi(a) +t )   & q \\ 
 & \mathcal{Q_F}: & \hat{q} & q \\ 
  \underline{n=6} & \mathcal{H}:& N  & \id  \\ 
  \underline{n= 8}& \mathcal{O}:  & \id  & (u,v) \mapsto R(u,v) := v^{\sigma+2} +uv  +u^\sigma  \end{array}
\]
\end{center}
For the root group labeling $(u,\Theta,\theta)$ in the root group labeling with epimorphism data $(u,\Theta,\theta,K, \Lambda, \nu, k)$ one easily extend this map to a map $\eta_s$ from $u(s)$ to the underlying field $K$. Because this root group labeling was in a standard form, this can be done unambiguously. 

The last element of the root group labeling with epimorphism data consists of a function $k:S \to \Lambda$, mapping an involution $s$ in $S$ to some element in $\nu(\eta_s(u(s)))$ different from $\infty$.

\end{itemize}

\subsection{Known results on epimorphisms of spherical buildings}\label{section:known}
Epimorphisms of generalized $n$-gons are well studied for generalized triangles (also known as \emph{projective planes}). Skornyakov expressed in~\cite{Sko:57} epimorphisms in terms of the coordinatizing planar ternary rings as \emph{places}. Subsequently the epimorphisms of projective Moufang planes and spaces have been classified (see~\cite{And:69},~\cite{Fau-Fer:83} and~\cite{Kli:56}). 

For other generalized polygons much less is known. There is a result of Pasini (\cite{Pas:84}) which says that the cardinalities of the preimages of an epimorphism between generalized $n$-gons are either always 1 or always infinite. This implies that epimorphisms between finite generalized $n$-gons are always isomorphisms. Epimorphisms from a generalized $n$-gon to a generalized $m$-gon with $m < n$ are studied by Gramlich and Van Maldeghem in~\cite{Gra-Mal:00} and~\cite{Gra-Mal:01}.

For other Moufang spherical buildings the only result known to the author are constructions using the theory of affine buildings and their non-discrete generalizations $\R$-buildings (see~\cite{Par:00} and~\cite{Wei:09}). One spherical building is then the `building at infinity' of an $\R$-building and the other a residue of it. We will call such morphisms \emph{affine epimorphisms}. The $\R$-buildings with an irreducible Moufang spherical building of rank at least 2 at infinity have been classified by Tits (see~\cite{Bru-Tit:72} and~\cite{Tit:86}). Without going in details, $\R$-buildings arise from valuations of the underlying (alternative) division algebra.

\begin{rem} \rm
The trivial epimorphisms, i.e. isomorphisms, can and will be considered to be affine epimorphisms in this paper. See Section~\ref{section:zero}.
\end{rem}

\begin{rem} \rm
In the non-Moufang case a wild variety of epimorphisms is possible. One way to do this is by using free constructions. Another way is to slightly perturbate the constructions of $\R$-buildings in~\cite{Str-Mal:11}, giving rise to epimorphisms of translation planes which are not arising from $\R$-buildings.
\end{rem}

\section{Statement of the main results and corollaries}\label{section:main}
The first Main Result shows that being Moufang is preserved under epimorphisms (note that this is trivial for higher dimensions).

\begin{Theorem}
The epimorphic image of a Moufang generalized polygon is again a Moufang polygon.
\end{Theorem}

The second and third Main Result classify the epimorphisms of a large class of spherical Moufang buildings.

\begin{Theorem}\label{Theorem:second}
Let $\Delta$ be spherical Moufang building of type listed in Section~\ref{section:deffield} with underlying field $K$. If $\phi : \Delta \to \Delta'$ is an epimorphism of buildings, then this gives rise to a root group labeling with epimorphism data $(u,\Theta,\theta,K, \Lambda, \nu, k)$, as defined in Section~\ref{section:rgl}, satisfying the compatibility conditions listed in Section~\ref{section:cc}. 
\end{Theorem}

\begin{Theorem}\label{Theorem:third}
If $(u,\Theta,\theta,K, \Lambda, \nu, k)$ is a root group labeling with epimorphism data  satisfying the compatibility conditions listed in Section~\ref{section:cc}, then one construct an epimorphism of buildings $\varphi:\Pi \to \Pi'$ where $\Pi$ is the Mou\-fang building associated to the root group labeling $(u,\Theta,\theta)$. If the root group labeling with epimorphism data arose from an epimorphism $\phi: \Delta \to \Delta'$ as in Main Result 2, then there exists isomorphisms $\psi: \Delta \to \Pi$ and $\psi': \Delta' \to \Pi'$ such that the following diagram commutes.
$$\xymatrix{ \Delta \ar[r]^{\psi}  \ar[d]_\phi &  \Pi  \ar[d]^\varphi  \\
\Delta' \ar[r]_{\psi'}   &  \Pi'}$$
\end{Theorem}

The following corollaries indicate that the `primitive' epimorphisms are the affine ones. 

\begin{Cor}
If moreover the valuation $\nu$ in this root group labeling with epimorphism data has finite rank (which is always the case if the underlying field has finite transcendency degree), then the associated epimorphism can be realized by combining a finite number of affine epimorphisms. 
\end{Cor}

\begin{Cor}
If an epimorphism of a spherical Moufang building of type listed in Section~\ref{section:deffield}  is not decomposable in two proper epimorphisms (i.e. not isomorphisms), then it is an affine epimorphism.
\end{Cor}




\section{Reduction to the generalized polygon case}
The aim of this section is to show how one can obtain epimorphisms between generalized polygons from higher rank spherical buildings. This will turn out to be useful when studying Moufang buildings via their rank 2 residues.
Let $\phi$ be an epimorphism between spherical buildings $(\cC,\delta)$ and $(\cC',\delta')$ of type $(W,S)$.

\begin{lemma}\label{lemma:find2}
If $C$ and $D$ are two chambers of $(\cC,\delta)$ such that $\delta'(C^\phi, D^\phi) = s \in S$, then there exists a chamber $E$ in $\cC$ such that $D^\phi = E^\phi$ and $\delta(C,E)=s$.
\end{lemma}
\proof
We start by finding a chamber $F$ of $(\cC,\delta)$ such that $F^\phi$ is opposite to both $C^\phi$ and $D^\phi$. As an epimorphism only can shorten the (numerical) distance between two chambers, one has that $F$ is opposite to both $C$ and $D$. If we project the chamber $D$ on the $s'$-panel containing $F$ (where $s'$ is the image of $s$ under the opposition involution) we obtain a chamber $G$ which is the unique chamber in this panel not opposite $D$. Clearly, its image is the unique chamber of the $s'$-panel containing $F^\phi$ not opposite to $D^\phi$. As $D^\phi$ is the unique chamber in the $s$-panel containing $C^\phi$ and $D^\phi$ not opposite to $G^\phi$, we have that the projection of the chamber $G$ back on the $s$-panel containing $C$ yields a chamber $E$ whose image has to be $D^\phi$. As $\delta(C,E)$ has to be $s$ by the definition of an $s$-panel, one has proven the lemma. \qed  


\begin{lemma}\label{lemma:redu}
Let $C$ be a chamber of $(C,\delta)$ and $S'\subset S$ a subset of size 2. Then $\phi$ induces an epimorphism from the $S'$-residue of $(C,\delta)$ containing $C$ to the $S'$-residue of $(C',\delta')$ containing $C^\phi$.
\end{lemma} 
\proof
The restriction of $\phi$ to the $S'$-residue of $(C,\delta)$ containing $C$  will map elements into the $S'$-residue of $(C',\delta')$ containing $C^\phi$ by the definition of epimorphisms and residues. Surjectivity of this morphism is a consequence of the previous lemma. \qed

\section{Proof of the first Main Result}\label{section:descend}

\subsection{Setting}\label{section:setting1}
Let $\Gamma:=(\cP,\cL,\inc)$ and $\Gamma'=(\cP',\cL',\inc')$ be two generalized $n$-gons, $\phi: \Gamma \to \Gamma'$ an epimorphism between them. Choose a root $\alpha:= (x_0, x_1, \dots , x_{n} )$ of $\Gamma$ which does not collapse under $\phi$. (To verify that these indeed exist pick $x_0$ and $x_n$ to be two elements of $\Gamma$ which are mapped to opposite elements, then each root beginning in $x_0$ and ending in $x_n$ cannot collapse as epimorphisms only shorten distances.) 

The main part of the proof is devoted to investigating under which conditions root elations of $\alpha$ descend under $\phi$.

Let $g$ be a root elation of $\alpha$. It maps an element $x_{-1}$ incident with $x_0$ but different from $x_1$ to an element $x_{-1}^g$ such that $x_{-1} \neq x_1 \neq x_{-1}^g$. We will prove that a sufficient condition for the root elation to descend is that there exists a neighbor $x_{-1}$ of $x_0$ such that $x_{-1}^\phi$  and $x_{-1}^{g \phi}$ are both different from $x_1^\phi$ (which is clearly a necessary condition as well). 

Once this is established, the first Main Result will follow quickly.

\subsection{Proof}
We start with an auxiliary lemma. 

\begin{lemma}\label{lemma:find}
If $a^\phi \inc' b^\phi$, then there exists an element $b'$ such that $b^\phi = b'^\phi$ and $a\inc b'$.
\end{lemma}
\proof
This is a reformulation of Lemma~\ref{lemma:find} in the language of epimorphisms between generalized polygons. \qed

We say that an element $x$ of $\Gamma$ has Property ($*$) if for each two elements $a, b \inc x$ one has that $a^\phi =b^\phi$ if and only if $a^{g \phi} =b^{g \phi}$.

\begin{prop}\label{prop:done}
If each element of $\Gamma$ has Property ($*$), then the root elation descends.
\end{prop}
\proof
First of all note that if two elements $a$ and $b$ of $\Gamma$ have opposite images under $\phi$, that then $a^g$ and $b^g$ also have opposite images under $\phi$ (because if a path does not collapse under $\phi$, then its image under $g$ will neither by Property ($*$)).

Suppose that the root elation does not descend, or equivalently that there exist elements $x$ and $y$ in $\Gamma$ such that $x^\phi =y^\phi$, but $x^{g\phi}  \neq y^{g \phi}$. (One also needs to consider the reverse statement, but this follows from an analogous exposition for the root elation $g^{-1}$.) We can assume that $x$ and $y$ are chosen so that the distance $k$ between them is minimal. Note that $k$ has to be bigger than zero and even, as $\phi$ does not map points to lines or vice versa. If $k$ would be 2, then $a$ and $b$ are both incident with some element $c$. Property ($*$) for this element then gives rise to a contradiction. 

Let $(y_0 := x,y_1, \dots ,y_k := y)$ be a path of shortest length between $x$ and $y$. Note that $x^\phi =  y_i^\phi$ only if $i=0$ or $k$, as otherwise it would contradict the way we choose the elements $x$ and $y$ (as it is impossible that $x^{g\phi} =y_i^{g \phi} =y^{g \phi}$). In particular this implies that $y_1^\phi= y_{k-1}^\phi$. Minimality of $k$ yields $y_1^{g\phi} =y_{k-1}^{g\phi}$.

Using Lemma~\ref{lemma:find}, one can find a path $(a_0, a_1, \dots, a_{n-1}:=x, a_n:=y_1)$ of length $n$ which does not collapse under $\phi$. So $a_0^\phi$ is opposite $y_1^\phi$. Combining this with $y_1^\phi =y_{k-1}^\phi$ gives that $a_0$ is opposite $y_{k-1}$. Let $(b_0:=a_0, b_1, \dots, b_{n-1}:=y, b_n:=y_{k-1}$) be the unique shortest path from $a_0$ to $y_{k-1}$ containing $y$ (which cannot collapse either). As $x^\phi$ equals $y^\phi$, and $y_1^\phi =  y_{k-1}^\phi$ is opposite to $a_0^\phi$, one has that $a_1^\phi = b_1^\phi$. As the distance between $a_1$ and $b_1$ is at most 2, this implies that $a_1^{g \phi} = b_1^{g \phi}$. Now because $y_1^{g\phi} =y_{k-1}^{g\phi}$ is opposite to $a_0^{g\phi}$, we have that the distance between $y_1^{g\phi}$ and $a_1^{g\phi}$ is $n-1$. In particular it follows that $x^{g\phi} =y^{g\phi}$, which contradicts the way we have chosen $x$ and $y$. \qed



\begin{lemma}\label{lemma:5.3}
Let  $(y_0, \dots, y_n)$ be a path of length $n$ in $\Gamma$ which does not collapse under $\phi$ and $\phi \circ g$. If Property ($*$) is satisfied for $y_n$, then it is also satisfied for $y_0$. 
\end{lemma}
\proof
Let $a$ and $b$ be two elements incident with $y_0$. Note that $y_0^\phi$ and $y_0^{g\phi}$ are opposite, respectively, to $y_n^\phi$ and $y_n^{g\phi}$. Because of this one has that $a^\phi = b^\phi$ if and only if $(\proj_{y_n}a)^\phi = (\proj_{y_n}b)^\phi$, and $a^{g\phi} = b^{g\phi}$ if and only if $(\proj_{y_n}a)^{g\phi} = (\proj_{y_n}b)^{g\phi}$. Property ($*$) for $y_n$ now implies that the conditions  $(\proj_{y_n}a)^\phi = (\proj_{y_n}b)^\phi$ and $(\proj_{y_n}a)^{g\phi} = (\proj_{y_n}b)^{g\phi}$ are equivalent. We conclude that $a^\phi = b^\phi$ if and only if $a^{g\phi} = b^{g\phi}$, so $y_0$ satisfies Property ($*$).\qed


\begin{cor}\label{cor:5.4}
If all elements of the root $\alpha$ satisfy Property ($*$), then all elements in every apartment which contains $\alpha$ and has an apartment as image, satisfy Property ($*$).
\end{cor}
\proof¤
Let $\Sigma$ be such an apartment. Then $\Sigma^{g\phi}$ is again an apartment by Property ($*$), and hence we can apply Lemma~\ref{lemma:5.3} to obtain that all elements of it satisfy this property.
\qed

Let $\Sigma$ be the unique apartment containing $x_{-1}$, $x_0$, \dots, $x_n$; and $\Sigma'$ the unique apartment containing $x_{-1}^g$, $x_0$, \dots, $x_n$. So $g$ maps $\Sigma$ to $\Sigma'$. Note that our assumption $x_{-1}^\phi  \neq x_1^\phi  \neq x_{-1}^{g\phi}$ implies that both apartments do not collapse under $\phi$. Let $x_{n+1}$ be the unique element of $\Sigma$ opposite $x_1$.

\begin{prop}\label{prop:5.5}
All the elements of $\Gamma$ satisfy Property ($*$).
\end{prop}
\proof
The elements $x_1, \dots, x_{n-1}$ all satisfy Property ($*$) as all elements incident with one of them are fixed by $g$. Applying  Lemma~\ref{lemma:5.3}, one then has that all elements of $\Sigma$, except from possibly $x_0$ and $x_n$, satisfy Property ($*$). Using Lemma~\ref{lemma:find} one can find an element $y_2 \inc x_1$ such that $x_0^\phi \neq y_2^\phi \neq x_2^\phi$. Let $(x_1,y_2, y_3, \dots, y_n, x_{n+1})$ be the unique shortest path from $x_1$ to $x_{n+1}$ containing $y_2$. Note that the path obtained by adding $x_0$ or $x_2$ as first element cannot collapse under $\phi$ or $\phi \circ g$ by the oppositeness of $x_1^\phi$ and $x_{n+1}^\phi$, and Property ($*$) for $x_1$. Lemma~\ref{lemma:5.3} applied to the path $(y_n, y_{n-1}, \dots, y_2, x_1, x_2)$ implies that $y_n$ satisfies Property ($*$), and applied to the path $(x_0,x_1,y_2,  \dots, y_n)$ it implies that $x_0$ satisfies Property ($*$). One concludes that all elements of $\Sigma$ satisfy Property ($*$).

Choose an element $z$ of $\Gamma$. Let $(z,z_1, \dots, z_k)$ be a shortest path from $z$ to an element $z_k$ of $\Sigma$ (`shortest' over all elements of $\Sigma$). There are exactly two apartments of $\Gamma$ containing a root of $\Sigma$ and the element $z_{k-1}$. As it is impossible that both apartments collapse under $\phi$ (this would imply that $\Sigma$ collapses as well), let $\Sigma''$ be such an apartment which does not collapse. It is easily seen that $(\Sigma'')^{g\phi}$ is again an apartment as $z_k$ satisfies Property ($*$). So by Corollary~\ref{cor:5.4} all elements of it satisfy property ($*$). By repeating this algorithm (substituting the role of $\Sigma$ by $\Sigma''$) a finite number of steps, one sees that $z$ satisfies Property ($*$). Hence all elements of $\Gamma$ satisfy Property ($*$). \qed

\begin{cor}\label{cor:5.6}
Let $\Sigma$ be an apartment of $\Gamma$ such that $\Sigma^\phi$ is an apartment of $\Gamma'$, number the elements of $\Gamma$ and its corresponding root groups in the usual way as described in Section~\ref{section:Moufang}. Then the following hold:
\begin{enumerate}[(a)]
\item $g\in U_0$ descends if and only if $x_{-1}^{g\phi} \neq x_1^\phi$.
\item If $g$ descends, then the image of $g$ in $\Aut (\Gamma')$ is nontrivial if and only if $x_{-1}^{g\phi} \neq x_{-1}^\phi$.
\end{enumerate}
\end{cor}
\proof
Part (a) follows by combining Propositions~\ref{prop:done} and~\ref{prop:5.5}. The `if' of part (b) is trivial, the `only if' follows from~\cite[Thm. 4.4.2 (v), (vi)]{Mal:98} as the image of $g$ then fixes $\Sigma^\phi$ and every element incident with $x_1^\phi$ or $x_2^\phi$.\qed

The first Main Result now follows easily.

\begin{cor}
The epimorphic image of a Moufang polygon is again a Mou\-fang polygon.  
\end{cor}
\proof
For every root $\alpha'$ in $\Gamma'$ one can find a root $\alpha$ in $\Gamma$ mapped to it using Lemma~\ref{lemma:find}. Even stronger, one can find for each two apartments $\Xi$ and $\Xi'$ containing $\alpha'$ two corresponding apartments $\Sigma$ and $\Sigma'$ in $\Gamma$. The unique root elation mapping $\Sigma$ to $\Sigma'$ descends as it has to satisfy the condition stated in Section~\ref{section:setting1}. Hence there is a root elation of $\alpha'$ mapping $\Xi$ to $\Xi'$. We conclude that $\Gamma'$ is a Moufang polygon.
\qed

\section{Epimorphisms and root groups}\label{section:rgd}
In this section we study various general properties that the root groups of a generalized Moufang polygon with an epimorphism should have. The main goal is to develop tools to be used in the next section where we invoke the classification of Moufang polygons and separate into cases.

\subsection{Subgroup information}\label{section:subg}
Let $\Sigma$ be an apartment of a generalized Moufang $n$-gon $\Gamma$ which does not collapse under an epimorphism $\phi: \Gamma \rightarrow \Gamma'$, and label the elements of it by $x_i$, with $i \in \mathbb{Z}$ such that $x_i \inc x_{i+1}$ and $x_i = x_{i+2n}$. This apartment will be called the \emph{hat-rack}. Let $U_i$ be the root group of the root $(x_i, x_{i+1}, \dots , x_{i+n})$. All of the $U_i$ forms the \emph{root group data} of $\Gamma$ associated to $\Sigma$. We will often use subscripts to indicate to which root group an automorphism belongs.

By Section~\ref{section:descend}, we have for each root group $U_i$ two subgroups $W_i \lhd V_i \leq U_i$ such that $V_i$ consists of all root elations of $U_i$ that descend, and $W_i$ consists of those that descend to the trivial automorphism. This implies that the root groups of the apartment $\Sigma^\phi$ of $\Gamma'$ are the quotients $U'_i := V_i / W_i$.

Later on, in Section~\ref{section:rigid} we will see that knowledge of the subgroups $W_1 \lhd V_1 \leq U_1$ and $W_n\lhd V_n \leq U_n$ suffices to determine the epimorphism uniquely. 

These results can be used to study a higher rank building $\Delta$ as well, this as the rank 2 residues containing a chamber $C$ in $\Delta$ form generalized polygons, on which epimorphisms are induced (see Lemma~\ref{lemma:find2}). Considering how the root group labeling $(u,\Theta,\theta)$ of $\Delta$ arises (see Section~\ref{section:high}), in particular the involvement of root group sequences of the rank 2 residues containing $C$, this implies the existence of subgroups  $w(s) \lhd v(s) \leq u(s)$. In Section~\ref{section:rigid} we will again show that this  suffices to determine the epimorphism uniquely.

\subsection{Opposite root groups }\label{section:opp}

In this section we investigate the behavior of two opposite root groups $U_i$ and $U_j$ in $\Sigma$ (meaning that $j\equiv i+n \mod 2n$). Without loss of generality we can assume that these are the root groups $U_0$ and $U_{n}$, that both fix the element $x_0$. Especially we consider the action of them on the elements incident with $x_0$. 

\begin{rem} \rm
This kind of action is also known as a Moufang set, for a detailed discussion see~\cite{Med-Seg:09}. 
\end{rem}

For an element $g \in U_n^*$ we define $\kappa_n(g)$ to be the unique element of $U_0$ which maps $x_{-1}$ to $x_1^g$. This defines a bijection from $U_n^*$ to $U_0^*$. 


\begin{lemma}\label{lemma:kappa}
The bijection $\kappa_n$ maps (bijectively)
\begin{itemize}
\item
$W_n^* $ to $U_0 \setminus V_0$, 
\item
$V_n \setminus W_n$ to $V_0 \setminus W_0$, 
\item
$U_n \setminus V_n$ to $W_0^*$.
\end{itemize}
\end{lemma}
\proof
Let $y = x_{-1}^{\kappa_n (g)\phi} =x_1^{g\phi}$. By repeated applications of Corollary~\ref{cor:5.6}, we see that if $g \in W_n^*$, then $y=x_1^\phi$ and hence $\kappa_n(g) \notin V_0$; if $g \in V_n \setminus W_n$, then $y \notin \{x_1^\phi, x_{-1}^\phi\}$ and hence $\kappa_n(g) \in V_0 \setminus W_0$; and if $g \notin V_n$, then $y = x_{-1}^\phi$, so $y \neq x_1^\phi$ because $\Sigma^\phi$ is an apartment, and hence $\kappa_n(g) \in W_0^*$. \qed

%
%
%


\begin{lemma}\label{lemma:switch}
Let $v_n$ be an element of $V_n \setminus W_n$. Then the map $g \in U_n \setminus V_n \mapsto {\kappa_n^{-1}(\kappa_n(v_n) \kappa_n(g))}v^{-1}_n$ is a bijection from $U_n \setminus V_n$ to $W_n^*$.
\end{lemma}
\proof
The orbit of $x_1^{v_n}$ under $W_0$ or $W_n$ is the preimage of $x^{v_n \phi}$ under $\phi$. In particular these orbits coincide. Also note that the groups $W_0$ and $W_1$ act regularly on the orbit. So we can conclude that we have a bijection which maps a $w_0 \in W_0$ to the unique element $w_n \in W_n$ such that $x_1^{v_n w_0} = x_1^{w_n v_n}$ (where we made use of the fact that $W_n \lhd V_n$). By the definition of $\kappa_n$ we now have that 
\begin{align*}
x_{1}^{w_n v_n} &= x_1^{v_n w_0 }   \\
 & = x_{-1}^{\kappa_n(v_n) w_0 } \\
 &= x_1^{\kappa_n^{-1}(\kappa_n(v_n) w_0)} .
\end{align*}
Hence $w_n = {\kappa_n^{-1}(\kappa_n(v_n) w_0)}v^{-1}_n$. The lemma is now proven because $\kappa_n$ is a bijection from  $U_n \setminus V_n$ to $W_0^*$ by the previous lemma. 
\qed

\subsection{Other pairs of root groups}

We now investigate the behavior of non-opposite root groups of the hat-rack. In particular we want to study the interacting with the commutator relations between them (see Section~\ref{section:Moufang}).

Recall that $U_{[i,j]}$ is the group generated by $U_i, U_{i+1}, \dots ,U_j$ if $i \leq j$ and the trivial group otherwise. We use similar notations $V_{[i,j]}$ and $W_{[i,j]}$ to denote the subgroup of generated by the subgroups of the form $V_k$ and $W_k$ respectively.


\begin{lemma} \label{lemma:prod}
If $i \leq j \leq i+n -1$, then the product $u_i u_{i+1} \dots u_j$ (where $u_k \in U_k$) descends if and only if each of the factors descend.
\end{lemma}
\proof
We prove this by induction. Assume that the product $g := u_i u_{i+1} \dots u_j$ descends. If $i =j$, then it is trivial that the factors descend, so suppose $i < j$. Note that $x_{j-1}^g = x_{j-1}^{u_j}$, hence $x_{j-1}^{u_j \phi} = x_{j-1}^{g \phi}  \neq x_{j+1}$ (the inequality holds as $g$ descends and fixes $x_{j+1}$). Corollary~\ref{cor:5.6} implies that $u_j$ descends. The product $g u_j^{-1} = u_i \dots u_{j-1}$ descends as well, so by induction all factors descend. The other direction is trivial.
\qed 

\begin{cor}\label{cor:prod}
If $i+1 \leq j \leq i+n-1$, then 
\begin{align*}
[V_i, V_j] & \leq V_{[i+1,j-1]}, \\ 
[V_i, W_j] & \leq W_{[i+1,j-1]}, \\
[W_i, V_j] & \leq W_{[i+1,j-1]}. 
\end{align*}
\end{cor}
\proof
Let $u_i \in U_i$, and $u_j \in U_j$. By  Lemma~\ref{lemma:nil2}, one can write $[u_i,u_j]$ in a unique way as a product $u_{i+1} u_{i+2} \dots u_{j-1}$, with $u_k \in U_k$. Now suppose that $u_i \in V_i$ and $u_j \in V_j$, then the product $u_{i+1} u_{i+2} \dots u_{j-1}$ descends, so Lemma~\ref{lemma:prod} implies that $u_{i+1} u_{i+2} \dots u_{j-1} \in V_{[i+1,j-1]} $. If moreover either $u_i \in W_i$ or $u_j \in W_j$, then their commutator descends to the trivial automorphism of $\Gamma'$. So the product $u_{i+1} \dots u_{j-1}$ is an element of $W_{[i+1,j-1]}$ by applying Lemma~\ref{lemma:nil2} to $\Gamma'$.
\qed

\subsection{Action of $\mu$-maps}
The $\mu$-maps form another type of interaction between the root groups, as the next lemma describes.

\begin{lemma}[\cite{Tit-Wei:02}, Prop. 6.1-2]
Let $\kappa_i : U_i^* \rightarrow U_{i+n}^*$ be as in Section~\ref{section:opp}. The automorphism $\mu_i(u_i) := \kappa_i(u_i) u_i (\kappa(u_i^{-1}))^{-1}$ (with $u_i \in U_i^*$) fixes $x_i$ and $x_{i+n}$, reflects $\Sigma$, and $U_j^{\mu_i(u_i)} = U_{2i+n -j}$ for each $j \in \mathbb{Z}$. \qed
\end{lemma}

Applying Lemma~\ref{lemma:kappa} this yields the following direct corollary.

\begin{cor}\label{cor:switch}
Let $v_i \in V_i\setminus W_i$, then 
\begin{align*}
V_j^{\mu_i(v_i)} &= V_{2i+n -j}, \\
W_j^{\mu_i(v_i)} &= W_{2i+n -j}
\end{align*}
for each $j \in \mathbb{Z}$. \qed
\end{cor}

The action of various $\mu$-maps can be found explicitly in~\cite[\S 32]{Tit-Wei:02}, and implicitly using~\cite[Lem. 6.4]{Tit-Wei:02}.

\begin{lemma}\label{lemma:prod2}
Choose a $u_1 \in U_1$ and a $u_n \in V_{n} \setminus W_{n}$.  Let $[u_1,u_n^{-1}] = u_2 \dots u_{n-1}$ (with $u_i \in U_i$), then 
\begin{align*}
u_1 \in V_1 &\Leftrightarrow  u_2 \in V_2, \\
u_1 \in W_1 &\Leftrightarrow  u_2 \in W_2.
\end{align*}
\end{lemma} 
\proof
Corollary~\ref{cor:prod} states that the implications from left to right are true. So suppose that $u_1 \in U_1 \setminus V_1$. By~\cite[Lem. 6.2, 6.4]{Tit-Wei:02} one has that $[u_2, \kappa_1(u_1^{-1})] = u_3 \dots u_{n-1} u_n$. As $U_1 \setminus V_1$ is stabilized under inversion as $V_1$ is a subgroup, it follows by Lemma~\ref{lemma:kappa} that $\kappa_1(u_1^{-1})\in W_{n+1}$. Using Corollary~\ref{cor:prod} and the assumption that $u_n \in V_{n} \setminus W_{n}$ yields that $u_2 \in U_{2} \setminus V_{2}$. This proves $u_1 \in V_1 \Leftrightarrow  u_2 \in V_2$. The proof of the second part is analogous.
\qed

\subsection{Rigidity and factorizations}\label{section:rigid}
We end this section by stating results on how the epimorphism is determined when certain $V_k$ and $W_k$ are known, and how different epimorphisms are related.


\begin{lemma}\label{lemma:rigid}
Let $\omega := u_2 \dots u_n$ with $u_i \in U_i$ for $i \in \{2, \dots, n\}$. The image of the element $x_1^\omega$ under the epimorphism $\phi$ is opposite $x_{n+1}^\phi$ if and only if all the root elations $u_i$ descend.
\end{lemma}
\proof
If each of the factors $u_i \in U_i$ ($i \in \{2, \dots, n\}$) descend then the product descends as well, so $(x_1^\omega)^\phi$ will be opposite $(x_{n+1}^\omega)^\phi = x_{n+1}^\phi$. Now suppose that one of the factors does not descend and let $u_j$ be the one with maximal index $j$, then one has by Corollary~\ref{cor:5.6} that $(x_{j-1}^\omega)^\phi = (x_{j-1}^{u_j \cdots u_n})^\phi = (x_{j+1}^{ u_{j+1} \cdots u_n})^\phi=(x_{j+1}^\omega)^\phi$, so the path $(x_1^\omega,x_2^\omega, \dots, x_{n+1}^\omega)$ collapses, or equivalently $(x_1^\omega)^\phi$ is not opposite $x_{n+1}^\phi$. \qed


\begin{prop}\label{prop:factor}
Suppose we have two different families of subgroups (given respectively by subgroups $W_k \lhd V_k < U_k$ and  $W'_k \lhd V'_k < U_k$, with $k=1$ or $n$) arising from two epimorphisms $\phi:\Gamma \rightarrow \Gamma_1$ and $\phi' :\Gamma \rightarrow \Gamma_2$  of the generalized $n$-gon $\Gamma$ with respect to the same hat-rack. If $V_k' \leq V_k$ and $W_k \leq W_k'$ for $k = 1,n$, then there exists an epimorphism $\phi'': \Gamma_1 \rightarrow \Gamma_2$ such that $\phi' = \phi'' \circ \phi$.
\end{prop}
\proof
By repeated applications of Corollary~\ref{cor:switch} (using elements of $V'_k \setminus W'_k$ for an appropriate index $k$), similar inclusions hold for the other root groups corresponding with the hat-rack. 

First we show that for each two elements $y$ and $z$ of $\Gamma$ one has that $y^\phi = z^\phi $ implies $ y^{\phi'} = z^{\phi'}$. The hat-rack in $\Gamma$ contains at least one element $x$ such that $x^{\phi'}$ is opposite to $y^{\phi'}$. Without loss of generality we may assume that this is the element $x_{n+1}$. In particular this element $x_{n+1}$ is opposite to $y$. Let $\omega := u_2 \dots u_n$ and $\omega' := u'_2 \dots u'_n$ (written as products of root elations with $u_i,u_i' \in U_i$) be the unique elements in $U_{[2,n]}$ such that $y = x_1^{\omega}$ and $z=x_1^{\omega'}$.  Note that we can do this because of Lemma~\ref{lemma:nil2} and as the group $U_{[2,n]}$ fixes the elements $x_n$ and $x_{n+1}$ of the hat-rack while acting regularly on elements opposite $x_{n+1}$. By Lemma~\ref{lemma:rigid} we have that $u_i \in V'_i$ as $y^{\phi'}$ is opposite to $x_{n+1}^{\phi'}$. As $V_k' \subset V_k$, it follows that $y^\phi$ is opposite to $x_{n+1}^{\phi}$. Because $y^\phi = z^\phi$, we obtain that $ \omega (\omega')^{-1} \in W_{[2, n]}$. Using $W_i \leq W'_i$,  we conclude that $y^{\phi'} = z^{\phi'}$. 

This property enables us to construct a surjective map $\phi''$ such that $\phi' = \phi'' \circ \phi$. The only thing left to prove is that $\phi''$ preserves adjacency. Let $a_1$ and $b_1$ be two incident elements of $\Gamma_1$. By Lemma~\ref{lemma:find} there exist incident elements $a$ and $b$ of $\Gamma$ such that $a_1= a^\phi$ and $b_1 = b^\phi$. It follows that $a_1^{\phi''} = a^\phi$ and $b_1^{\phi''} = b^\phi$ are incident.
\qed


\begin{cor}\label{cor:rigid}
If the subgroups $V_1$, $V_n$, $W_1$ and $W_n$ are known, then the epimorphism $\phi$ is unique (up to an isomorphism of the image).
\end{cor}
\proof
Proposition~\ref{prop:factor} implies that if there are two epimorphisms $\phi:\Gamma \rightarrow  \Gamma_1$ and $\phi': \Gamma \rightarrow \Gamma_2$ with the same subgroups $V_1$, $V_n$, $W_1$ and $W_n$ (with respect to the same hat-rack), that then there exist epimorphisms $\phi'' : \Gamma_1 \rightarrow \Gamma_2$, $\phi''':  \Gamma_2 \rightarrow \Gamma_1$ such that $\phi' = \phi'' \circ \phi$ and $\phi = \phi''' \circ \phi'$. One easily verifies that $\phi''$ and $\phi'''$ are inverses of each other, hence we obtain that the epimorphism is unique up to an isomorphism.
\qed

\begin{cor}\label{cor:rigid2}
If the subgroups $w(s) \lhd v(s )\leq u(s)$ are known for a root group labeling $(u,\Theta,\theta)$ of a spherical Moufang building $\Delta$ (see Section~\ref{section:subg}), then the corresponding epimorphism of $\Delta$ is unique (up to an isomorphism of the image).
\end{cor}
\proof
This follows from the previous corollary and Tits' extension theorem~\cite[Thm. 4.2.1]{Tit:74}.
\qed

\section{Proof of the second Main Result}
In contrast with the previous section we now invoke the classification of irreducible spherical Moufang buildings of rank at least 2 and study what the properties determined in the previous section imply. This will lead us to a classification of epimorphisms of those buildings listed in Section~\ref{section:deffield}.

\textbf{Sketch of proof.} |
We start with assuming the existence of an epomorphism $\phi$. For deriving necessary conditions we look at the epimorphisms induced on the rank 2 residues (see Lemma~\ref{lemma:redu}).
 In Sections~\ref{section:psl} up to~\ref{section:ortho}, we study one pair of opposite root groups and show that the subgroups $V_i$ and $W_i$ arise from a valuation to an ordered abelian group of the underlying field. We then use this information to study the other root groups (Section~\ref{section:imp}) and determine certain conditions that need to be satisfied (Section~\ref{section:cc}). In Section~\ref{section:concl2} we use this information to construct the root group labeling with epimorphism data and conclude the proof of the second Main Result.


\begin{rem} \rm
From this point on we only work with the root groups, not with elements of generalized polygons. In particular notations of the form $x_i$ will now denote parametrizations of the root groups, not elements of a hat-rack as in the previous section. These parametrizations are of the form $x_i: M \to U_i$, where $M$ is some additive algebraic structure.
\end{rem}

\subsection{Projective lines} \label{section:psl}
In this section we assume that $U_0$ and $U_n$ are isomorphic to the additive group of an alternative division ring $K$ (later on we will restrict to fields), by maps $x_i: K \to U_i$ ($i \in \{1,n\}$), and that the map $\kappa_0$ is given by $x_0(a) \mapsto x_n(a^{-1})$.

Applying Lemma~\ref{lemma:switch} we obtain that the map $\phi_a: x_0(b) \mapsto x_0((a^{-1} + b^{-1})^{-1} - a)$ is an bijection from $U_0 \setminus V_0$ to $W_0^*$ for every $x_0(a) \in V_0 \setminus W_0$. The expression $(x_0(b-a)^{\phi_a})^{-1}$ simplifies to $x_0(a b^{-1} a)$, which is also an involutory bijection from $U_0 \setminus V_0$ to $W_0^*$, as $V_0$ and $W_0$ are subgroups of $U_0$. 

Choose an element $t\in K$ such that  $x_0 (t) \in V_0 \setminus W_0$. Define the following subsets of $K$:
\begin{align*}
A &:= \{yt^{-1} \in K \vert x_0(y) \in U_0 \setminus V_0   \}, \\
B &:= \{yt^{-1} \in K \vert x_0(y) \in V_0 \setminus W_0   \}, \\
C &:= \{yt^{-1} \in K \vert x_0(y) \in W_0^* \}.
\end{align*}
\begin{lemma}
The subset $R := B \cup C \cup \{0\}$ forms a subring of $K$ containing the identity element.
\end{lemma}
\proof
Observe that $y \mapsto b y^{-1} b$ interchanges $A$ and $C$ bijectively for every $b \in B$. This implies that this map stabilizes $B$. As $1 = t\cdot t^{-1} \in B$, one also has that the inverse is a map of this form and that squaring stabilizes $B$. Lastly we note that $y \mapsto b y b$ stabilizes all three subsets $A$, $B$ and $C$  for every $b \in B$ (by combining the maps $y \mapsto b y^{-1} b$ and $y \mapsto y^{-1}$).

The subset $R$ is an additive subgroup of $K$ as $V_0$ is a subgroup of $U_0$. So in order to show that $R$ forms a subring, we only need to show that it is closed under multiplication. First suppose that $b$ and $c$ both lie in $B$. The maps $y \mapsto b^{-1} y  b^{-1} $ and  $y \mapsto c^{-1} y c^{-1}$ stabilize the sets $A$, $B$ and $C$. The combination of both maps $bc$ to $(cb)^{-1}$. If $b$ and $c$ commute this implies that $bc \in B$ (as taking the inverse interchanges $A$ and $C$).  If $b$ and $c$ do not commute and $bc \notin B$ then the sum $bc  + cb$ lies in $A$ (as $bc$ and $cb$ cannot both lie in $A$ or $C$ at the same time). If $b+c \in B$, then we know that the square $(b+c)^2 =b^2 +c^2 +bc+cb$ also lies in $B$, but this is a contradiction as $b^2$ and $c^2$ are elements of $B$ while $bc+cb \in A$. If $b+c \in C$, then one can obtain a contradiction in a similar way considering $(1+b+c)^2$. We conclude that $bc \in B$. 

Now suppose that $b \in B$ and $c \in C$. So $1+c \in B$, hence by the previous paragraph we have that $b(1+c) =b + bc \in B \subset R$. As $R$ is closed additively, we have that $bc \in R$. 

The last case is handled analogously. Suppose that $b,c \in C$, then $1+c \in B$, so $b(1+c) =b +bc \in R$, hence again $bc \in R$.
\qed

\begin{lemma}
The set of units of $R$ is $B$, and $K = R \cup (R^*)^{-1}$.
\end{lemma}
\proof 
In order to prove this notice that taking the inverse stabilizes $B$, and interchanges $A$ and $C$.     
\qed

\begin{rem} \em
A subring $R$ of $K$ such that $K= R \cup (R^*)^{-1}$ is also known as a \emph{total subring}.
\end{rem}

\begin{cor}\label{cor:valu}
If $k$ is a field, then there exists a valuation $\nu$ of $K$ to an ordered abelian group $\Lambda$ and the symbol $\infty$ such that 
\begin{align*}
A &= \{y \in K \vert  \nu(y) < 0  \}, \\
B &= \{y \in K \vert  \nu(y) = 0 \}, \\
C &= \{y \in K \vert \nu(y) > 0 \}.
\end{align*}
\end{cor}
\proof
The previous lemma implies that $R$ is a valuation ring, and hence defines a valuation with the desired properties (see~\cite{Kru:32}). 
\qed

Returning to the root group $U_0$, we now have in the case that $K$ is a field that 
\begin{align*}
V_0 &= \{x_0(a) \in U_0 \vert  \nu(a) \geq l \}, \\
W_0 &= \{x_0(a) \in U_0 \vert  \nu(a) > l \},
\end{align*}
where $l = \nu(t)$. Recall that $\kappa_0$ maps $x_0(a)$ to $x_n(a^{-1})$. So Lemma~\ref{lemma:kappa} yields 
\begin{align*}
U_n \setminus W_n &= \{x_n(a) \in U_n \vert  \nu(a) \leq -l \}, \\
U_n \setminus V_n &= \{ x_n(a) \in U_n \vert  \nu(a) > -l \},
\end{align*}
which is equivalent to
\begin{align*}
V_n &= \{x_n(a) \in U_n \vert  \nu(a) \geq -l \}, \\
W_n &= \{ x_n(a) \in U_n \vert  \nu(a) > -l \}.
\end{align*}
\begin{rem} \em
Corollary~\ref{cor:valu} is not true for skew fields or octonion algebras, as a total subring is not necessarily stabilized by inner automorphisms, which is necessary for obtaining a valuation.
\end{rem}


\subsection{Orthogonal Moufang sets}\label{section:ortho}
The only case where there are no opposite root groups of the form discussed in Section~\ref{section:psl} are the Moufang quadrangles of exceptional type and those of indifferent type (so $n=4$). We approach these by considering a full subquadrangle of quadratic form type (\emph{full} means that we do not have to restrict the root groups of even index). The epimorphism of the entire quadrangle implies one of the subquadrangle, but not necessarily to a thick generalized quadrangle. The `full' property assures us at least some thickness, and due to the fact that the epimorphism arises by restricting root groups, one still can consider subgroups $V_k$ and $W_k$ and apply the results from Section~\ref{section:rgd}.

Let us describe this subquadrangle. Let $K$ be a field, $L_0$ a vector space over $K$ equipped with an anisotropic quadratic form $q: L_0 \to K$. Let $f$ be the bilinear form associated to $q$. Let the root groups $U_0$, $U_2$ and $U_4$ be parametrized by the additive group of $L_0$ via isomorphims $x_0$, $x_2$ and $x_4$. The root groups $U_1$, $U_3$ and $U_5$ are parametrized by the additive group of the field $K$ via isomorphisms $x_1$, $x_3$ and $x_5$. The map $\kappa_0 : U_0 \to U_4$ is given by $x_0(u) \mapsto x_4(u/q(u))$.   Because the subquadrangle is full, we have that $V_k \neq W_k$ for $k$ even. This is however not guaranteed for those of odd index (and hence we cannot apply the results from Section~\ref{section:psl} directly). We also list the non-trivial commutator relations between the root groups $U_1$, $U_2$, $U_3$ and $U_4$ (see~\cite[16.3]{Tit-Wei:02}):
\begin{align*}
[x_2(a), x_4(b)^{-1}] &= x_3(f(a,b)), \\
[x_1(t), x_4(a)^{-1}] &= x_2(ta) x_3(t q(a)).
\end{align*}
The existence of such a subquadrangle (and with similar notations) of the Moufang quadrangles of exceptional type $\mathsf{E}_i$ ($i =6, 7, 8$)  and $\mathsf{F_4}$ follows from the description~\cite[16.6-7]{Tit-Wei:02}, for those of indifferent type the notations from~\cite[16.4]{Tit-Wei:02} and our notations are related by the following table.

\begin{center}\begin{tabular}{c|c} 
Our notations &  \cite[16.4]{Tit-Wei:02} \\\hline 
$K$ & $K^2$  \\   
$L_0$  & $L_0$ \\
$q$ & $x \mapsto x^2$
\end{tabular} 
\end{center}

Let $x_4(a)$ be an element of $V_4 \setminus W_4$, and $b$ an element of $L_0$, linear independent of $a$. Denote by $\widehat{L}_0$ the two-dimensional subspace of $L_0$ spanned by both $a$ and $b$. There exists a quadratic extension $F/K$ with norm $N$ and a $K$-linear bijection $\theta$ from $F$ to $\widehat{L}_0$ such that $\theta(1) = a$ and $N(s) = q(\theta(s))/q(a)$ for all $s \in F$ (see for example~\cite[\S 2.6]{Med-Hao-Kno-Mal:06}). For $i = 0, 2$ and $4$ this subspace corresponds to a subgroup $\widehat{U}_i$ of $U_i$, parametrized by the map $x_i \circ \theta: F \to \widehat{U}_i$.


If the field extension $F/K$ is separable, then we denote by $\sigma$ the Galois involution of the extension. If it is inseparable, then $\sigma$ will be the identity.

Invoking Section~\ref{section:psl} on $\widehat{U}_0$ and $\widehat{U}_4$, we obtain a valuation $\omega$ of $F$ such that (with $l := \omega(\theta^{-1}(a)$):
\begin{align*}
V_4 \cap \widehat{U}_4  &= \{ x_4(\theta(y)) \in \widehat{U}_4 \vert  \omega(y) \geq l \}, \\
W_4 \cap \widehat{U}_4 &= \{ x_4(\theta(y))  \in \widehat{U}_4 \vert  \omega(y) > l \}.
\end{align*}
Note that the restriction of $\omega$ to $K$ does not depend of the choice of $b$. Also observe that each one-dimensional subspace of $\widehat{L}_0$ contains elements which are mapped to elements of $V_4$ by $x_4 \circ \theta$ (and analogously for $V_0$ and $V_2$).   

We now claim that the automorphism $\sigma$ arising from the field extension leaves the valuation $\omega$ invariant. Suppose that this is not the case, so there exists a $w \in F$ such that $\omega(w) < \omega(w^\sigma)$. Note that the field extension $F/K$ must be separable and accordingly that the bilinear form $f$ restricted to $\widehat{L}_0$ is non-trivial. Combined with the observation on one-dimensional subspaces of $\widehat{L}_0$, Corollary~\ref{cor:prod} and the commutator relations between $U_2$ and $U_4$ this yields that the subgroup $V_3 < U_3$ is non-trivial. Corollary~\ref{cor:switch} then implies that the subgroup $V_1 < U_1$ is non-trivial as well.

By the commutator relations and Lemma~\ref{lemma:prod2} we have that whenever $x_1(t) \in V_1$, then $ \{ x_2(\theta(y)) \in \widehat{U}_2 \vert  \omega(y) \geq l + \omega(t) \} \subset V_2 \cap \widehat{U}_2$. A consequence of this is $t$ cannot have arbitrary small valuations unless $\omega$ is the trivial valuation, as this would imply that $ \widehat{U}_2 \subset V_2$ for every choice of $b$ (and hence $U_2 =V_2$). A similar thing is true for choices of $x_3(t) \in V_3$ by Corollary~\ref{cor:switch}.

Let $d := w^{-1+\sigma}$, so $d^\sigma = d^{-1}$ and $\omega(d) > 0$. Hence $x_4(\theta(d^{m}))$ will be an element of $V_4$ for high enough values of $m$.  So for these values of $m$ and an arbitrary $c\in F^*$ such that $x_2(\theta(c)) \in V_2$ one has that the following commutator is contained in $V_3$ by Corollary~\ref{cor:prod}:
\begin{align*}
[x_2(\theta(c)), x_4(\theta(d^m))^{-1}] &= x_3(f(\theta(c),\theta(d^m))) \\
&= x_3(q(\theta(c) +\theta(d^m))  -q(\theta(c)) - q(\theta(d^m))   ) \\
&= x_3(N(c+d^m) - N(c) - N(d^m)) \\
&= x_3(c^\sigma (d^m) + c (d^m)^\sigma) \\
&= x_3(c d^m + c d^{-m}) \\
&= x_3(c (d^m + d^{-m})) .
\end{align*}
The last factor has an arbitrary low valuation using arbitrary large $m$. This contradicts the earlier remark that one cannot choose a $x_3(t) \in V_3$ with $t\in K$ having arbitrary small valuations.





We conclude that  $\sigma$ leaves $\omega$ invariant, so for an element $x_4(\theta(f)) \in \widehat{U}_4$ we have that $\omega(q(\theta(f))) = \omega(N(f)) = 2\omega(f)$. As the valuation $\omega$ restricted to $K$ is independent of the choice of $b$, we finally obtain:
\begin{align*}
V_4  &= \{ x_4(v) \in U_4 \vert  \omega(q(v)) \geq 2l \}, \\
W_4 &= \{ x_4(v) \in U_4 \vert  \omega(q(v))  > 2l \}.
\end{align*}

\subsection{Implications on the root group sequence}\label{section:imp}

Assume we have an epimorphism  $\phi: \Gamma \rightarrow \Gamma'$ between Moufang polygons. We use the description of the root group sequence as found in~\cite[\S{}16 and \S{}32]{Tit-Wei:02} (parametrizing the root group $U_r$ by a map $x_r$, with $r \in \{1, \dots,n\}$). One can define a `norm' function on the algebraic structure defining the other root groups into the underlying field. We list the functions in question in the following table:
\begin{center}
\[ \begin{array}{lccc} 
& & \underline{i \mbox{ odd}} &  \underline{i \mbox{ even} }\\
\\
\underline{n=3} & \mathcal{T}: & \id  & \id\\ 
\underline{n=4} & \mathcal{Q_Q}: & \id & q \\ 
& \mathcal{Q_D}: & a \mapsto a & a \mapsto a^2 \\ 
 & \mathcal{Q_P}: & (a,t) \mapsto t  & \id \\
 & \mathcal{Q_E}: & (a,t) \mapsto q(\pi(a) +t )   & q \\ 
 & \mathcal{Q_F}: & \hat{q} & q \\ 
  \underline{n=6} & \mathcal{H}:& N  & \id  \\ 
  \underline{n= 8}& \mathcal{O}:  & \id  & (u,v) \mapsto R(u,v) := v^{\sigma+2} +uv  +u^\sigma  \end{array}
\]
\end{center}
We will denote the `norm' function on $U_j$ by a generic $\eta_j$ regardless of type.  Fix $i$ to be $1$ if $\Gamma$ is of type $\mathcal{T}$, $\mathcal{Q_Q}$ or $\mathcal{O}$ and let $i=n$ in every other case, this for the rest of this section. Also set $j$ to be $2$ when $i = n$, and $n-2$ when $i=1$. The importance of the norm functions is illustrated by the following lemma. 

\begin{lemma}\label{lemma:easyprods}
Let $u_1 \in U_1$ and $u_n \in U_n$ be two root elations. If one writes $[u_1,u_n^{-1}]$ as a product $ u_2 \dots u_{n-1}$ ($u_r \in U_r$), then $\eta_j(u_j) = \pm \eta_1(u_1) \eta_n(u_n)$.
\end{lemma}
\proof
By straightforward calculations using the commutator relations found in~\cite[\S 16]{Tit-Wei:02}.
\qed
%

By applying the case studies made in Sections~\ref{section:psl} and~\ref{section:ortho} to the explicit descriptions in \cite[\S{}16 and \S{}32]{Tit-Wei:02} one observes that for the generalized polygons defined over a field listed in Section~\ref{section:deffield}, that then there exists a valuation $\nu :  K \twoheadrightarrow \Lambda \cup \{ \infty \}$ and $l \in \Lambda$ such that
\begin{align*}
V_i &= \{x_i(a) \in U_i \vert  \nu(\eta_i(a)) \geq l \}, \\
W_i &= \{x_i(a) \in U_i \vert  \nu(\eta_i(a)) > l \}.
\end{align*}

Choose a $v_{n+1-i} \in V_{n+1-i} \setminus W_{n+1-i}$, and let $k := \nu(\eta_{n+1-i}(x_{n+1-i}^{-1}(v_{n+1-i})))$. One is now able to describe the groups $V_j$ and $W_j$, and subsequently $V_{n+1-i}$ and $W_{n+1-i}$. 

\begin{lemma}\label{lemma:above}
\begin{align*}
V_j &= \{x_j(a) \in U_j \vert  \nu(\eta_j(a)) \geq k + l \}, \\
W_j &= \{ x_j(a) \in U_j \vert  \nu(\eta_j(a)) > k + l \}.
\end{align*}
\end{lemma}
\proof
We will prove this under the assumption that $i = n$ (so $j=2$), the other case is symmetric. Let $u_2 \in V_2$, and $u_n := u_2^{(\mu_1 (v_1)^{-1})}$. Using~\cite[Lem. 6.4]{Tit-Wei:02} this implies that $[v_1,u_n^{-1}] = u_2 u_3 \dots u_{n-1}$ with $u_r \in U_r$ for $r \in \{3, \dots, n-1\}$. The previous lemma yields that $\nu(\eta_2(x_2^{-1}(u_2))) = \nu(\eta_1(x_1^{-1}(u_1))) + \nu(\eta_n(x_n^{-1}(u_n)))= k + \nu(\eta_n(x_n^{-1}(u_n)))$. The statement now follows from Corollary~\ref{cor:switch}.
\qed

\begin{cor}
\begin{align*}
V_{n+1-i} = \{x_{n+1-i}(a) \in U_{n+1-i} \vert  \nu(\eta_{n+1-i}(a)) \geq k \}, \\
W_{n+1-i} = \{x_{n+1-i}(a) \in U_{n+1-i} \vert  \nu(\eta_{n+1-i}(a)) > k \}. 
\end{align*}
\end{cor}
\proof
From Lemmas~\ref{lemma:prod2},~\ref{lemma:easyprods} and~\ref{lemma:above}. 
\qed

We now have a description of $V_1$, $V_n$, $W_1$ and $W_n$, which suffices to describe the epimorphism by Corollary~\ref{cor:rigid}.

The next goal is now to derive compatibility conditions from the commutator relations. We start by describing the other subgroups of interest of $U_r$ with $r \in \{2, \dots, n-1 \}$, using Corollary~\ref{cor:switch} (and the relations given in~\cite[\S{}16 and \S{}32]{Tit-Wei:02}) a finite number of times. We display this information schematically as a vector where the $n$ coordinates correspond to respectively $U_1, \dots ,U_n$, and the value at a coordinate $r$ is the element of $\Lambda$ which defines the subgroups $V_r$ and $W_r$ as a hyperlevel set and strict hyperlevel set respectively with respect to $\nu \circ \eta_r \circ x_r^{-1}$. 
\begin{center}
\[ \begin{array}{lcc}
 \underline{n=3} & \mathcal{T}: & (l, l+k ,k) \\
 \underline{n=4} & \mathcal{Q_P}:& (k,l+k, l + l'+k,l)\\
 & \mathcal{Q_Q}: & (l,2l+k,l+k, k) \\ 
 & \mathcal{Q_D}, \mathcal{Q_E}, \mathcal{Q_F}: & ( k ,l+k, 2l+k  ,l)  \\  
  \underline{n=6} & \mathcal{H}:& ( k,l+k, 3l+2k,2l+k, 3l+k,l)    \\ 
  \underline{n= 8} & \mathcal{O}: & (l,  2l+l' +k, l+l'+k,  2l+2l' +k+k' ,  \\ 
  & &  l+l'+2k-k',   2l+l' +k+k' ,l+k, k)
  \end{array}
\]
\end{center}

The $l'$ for the quadrangle of pseudo-quadratic form type $\mathcal{Q_P}$ and the octagon case $\mathcal{O}$ is defined as follows. Let $a \in K$ such that $\nu(a) = l$, then we set $l' := \nu(a^\sigma)$. Note that for these two cases we have $x_i(a) \in V_i \setminus W_i$ if and only if $\nu(a)=l$, and each such $x_i(a)$ can be used to apply Corollary~\ref{cor:switch} to. As the resulting $V_r$ and $W_r$ should be identical, the resulting $l'$ is independent of the $a \in K$ used to define it. The element $k'$ is defined in a similar way, so if $(u,v) \in K \times K$ such that $\nu(R(u,v)) = k$, then $k' := \nu((R(u,v))^\sigma)$.



\subsection{Compatibility conditions}\label{section:cc}
In this section we describe extra conditions who arise for epimorphisms of certain Moufang polygons.  These conditions, which we call the \emph{compatibility conditions}, involve the valuation $\nu$, the underlying algebraic structures and certain constants $k$ and $l$. We now list these compatibility conditions, the details for each case can be found in Sections~\ref{section:ccf}-\ref{section:ccl}. 
\begin{itemize}
\item Quadrangles $\mathcal{Q}_\mathcal{P}(K,K_0, \sigma, L_0,q)$ of pseudo-quadratic form type: 
\begin{align*}
&\forall t \in K: \nu(t) = \nu(t^\sigma), \\
&\forall (u,t),(v,s) \in T : \nu(t),\nu(s) \geq k \Rightarrow \nu(f(u,v)) \geq k.
\end{align*}
\item Quadrangles $\mathcal{Q}_\mathcal{Q}(K,L_0, q)$ of quadratic form type:
\[
\forall u,v \in L_0: \nu(q(u)),\nu(q(v)) \geq k \Rightarrow \nu(f(u,v)) \geq k.
\]
\item Quadrangles $\mathcal{Q}_\mathcal{E}(K,L_0, q)$, $\mathcal{Q}_\mathcal{F}(K,L_0, q)$ of types  $\mathsf{F}_4$, $\mathsf{E}_6$, $\mathsf{E}_7$ and $\mathsf{E}_8$: A list of ten conditions listed in Equation~\ref{eq:except} of Section~\ref{section:frank2}.
\item Hexagons $\mathcal{H}(J,F,\#)$:
\[
\forall u,v \in J : \nu(N(u)) \geq k,\nu(N(v)) \geq -k \Rightarrow \nu(T(u,v)) \geq 0.
\]

\item Octagons $\mathcal{O}(K,\sigma)$:
$$ \forall x \in K : \nu(x) = 0 \Rightarrow \nu(x^\sigma) = 0.$$
\end{itemize}

For an epimorphism of a Moufang building of type $(W,S)$ listed in Section~\ref{section:deffield}, these conditions occur if one considers the induced epimorphisms on the $\{s,t\}$-residues containing a certain chamber, where $s$ and $t$ is a pair of involutions in $S$. The role of the $k$ and $l$ are then replaced by constants $k(s)$ and $k(t)$, which will be defined in Section~\ref{section:concl2}. By the list of possible rank two residues given in Section~\ref{section:high} there is at most one choice of pair $\{s,t\}$ giving rise to compatibility conditions.

\begin{rem} \rm
We will not always derive the strongest conditions possible. This will not be a problem (and is even slightly beneficial) as we will see in Sections~\ref{section:oner} and~\ref{section:frank}.
\end{rem}

\subsubsection{Quadrangles $\mathcal{Q}_\mathcal{P}(K,K_0, \sigma, L_0,q)$ of pseudo-quadratic form type} \label{section:ccf}
The first compatibility condition involves the involution $\sigma$. From the appearance of the $l'$ in the last list, one has that if $x,y \in K$ and $\nu(x)=\nu(y)=l$, then $l'= \nu(x^\sigma) = \nu(y^\sigma)$. Note that there exists an $x\in K$ such that $\nu(x)=l$. Now suppose that $y,z$ are elements of $K$ such that $\nu(y) =\nu (z)$, then $\nu(xyz^{-1}) = l$. So $\nu((xyz^{-1})^\sigma) = \nu(x^\sigma)$, which implies that $\nu(y^\sigma) = \nu(z^\sigma)$. Suppose that there is an $y \in K$ such that $\nu(y) < \nu(y^\sigma)$, then $\nu(1+ y^{-1+\sigma}) = \nu(1)= 0$. Applying $\sigma$ yields 
\begin{align*}
0 = \nu(1^\sigma) &=\nu((1+ y^{-1+\sigma})^\sigma) \\
&= \nu(1+y^{1-\sigma}) \\
&= \nu(y^{1-\sigma}) <0,
\end{align*}
which is a contradiction. We conclude as first compatibility condition that $$\forall t \in K: \nu(t) = \nu(t^\sigma).$$ Note that this implies that $l' =l$.

A second compatibility condition involves the skew-hermitian form $f$. By the commutator relations between $U_1$ and $U_3$ (and Corollary~\ref{cor:prod}), we have that if $(u,t),(v,s) \in T$ with $\nu(t) \geq k$, $\nu(s) \geq 2l+k $. By substituting $(v,s)$ by $(vx,x^\sigma s x)$, where $\nu(x) =l$, one simplifies this to
\[
\forall (u,t),(v,s) \in T : \nu(t),\nu(s) \geq k \Rightarrow \nu(f(u,v)) \geq k.
\]
 
\subsubsection{Quadrangles $\mathcal{Q}_\mathcal{Q}(K,L_0, q)$ of quadratic form type}

In a similar way as for the second compatibility condition for pseudo-quadratic forms one obtains (using the commutator relations between $U_2$ and $U_4$) that
\[
\forall u,v \in L_0: \nu(q(u)),\nu(q(v)) \geq k \Rightarrow \nu(f(u,v)) \geq k.
\]

\begin{rem}\label{rem:comp}\rm
Let us consider the special case that the quadrangle is also of quadratic or honorary involutory type. So $L_0$ is a composition algebra equipped over $K$ equipped with (anisotropic) norm $q$. The map $\nu' := \nu \circ q$ on this composition division algebra satisfies $\nu'(u \cdot v) = \nu'(u)+ \nu'(v)$. We now claim that $\nu'$ is a valuation on $L_0$. The only thing we have still to verify is that if $\nu'(u) \leq \nu'(v)$ for $u,v \in L_0$, that then $\nu'(u) \leq \nu'(u+v)$. So let $u,v \in L_0$ with $\nu'(u) \leq \nu'(v)$. Choose an element $w \in L_0$ such that $\nu'(w)=k$. Observe that $\nu'(w\cdot v\cdot u^{-1})$ is at least $k$, so $\nu(f(w,w\cdot  v\cdot u^{-1})) \geq k$ by the above compatibility condition. Applying the relationship $f(u,v) =q(u+v) -q(u) -q(v) $ between a quadratic form and its associated bilinear form together with $\nu$ being a valuation yields that 
\begin{align*}
k &\leq \nu(q(w + w\cdot v\cdot u^{-1})) \\
&= \nu'(w+w\cdot v\cdot u^{-1}) \\
&= \nu'(w\cdot u^{-1}\cdot (u+v)) \\
&= \nu'(w) - \nu'(u) + \nu'(u+v) \\
&= k - \nu'(u) + \nu'(u+v),
\end{align*}
which simplifies to $\nu'(u) \leq \nu'(u+v)$, so $\nu'$ is indeed a valuation.
\end{rem}

\subsubsection{Quadrangles $\mathcal{Q}_\mathcal{E}(K,L_0, q)$, $\mathcal{Q}_\mathcal{F}(K,L_0, q)$ of types  $\mathsf{F}_4$, $\mathsf{E}_6$, $\mathsf{E}_7$ and $\mathsf{E}_8$} \label{section:ccx}

A list of the compatibility conditions (ten in total) one needs for these cases is listed in Equation~\ref{eq:except} in Section~\ref{section:frank2}, where $\phi_r$ is the function $\nu \circ \eta_r$ and $\eta_r$ as defined in Section~\ref{section:imp}. These follow from a combination of Corollary~\ref{cor:prod} and the results of Section~\ref{section:imp}.

The number of equations is much larger than in the other cases because the residues of the affine buildings associated to the generalized Moufang quadrangles of these types are not fully described yet. When such a description becomes available (as announced in~\cite[Rem. 21.43 and p. 228]{Wei:09}), one can expect to reduce the number of equations needed substantially.

%
%
%
%
%

\subsubsection{Hexagons $\mathcal{H}(J,F,\#)$}
Applying the same argument as for the second compatibility condition for pseudo-quadratic forms to the root groups $U_1$ and $U_3$ gives
\[
\forall u,v \in J : \nu(N(u)) \geq k,\nu(N(v)) \geq -k \Rightarrow \nu(T(u,v)) \geq 0.
\]

\subsubsection{Octagons $\mathcal{O}(K,\sigma)$} \label{section:ccl}
The sole compatibility condition for octogonal systems involves the Tits endomorphism $\sigma$. Analogously to the first compatibility condition for the pseudo-quadratic case one shows that if one has  two elements $x$, $y$ of $K$ such that $\nu(x)=\nu(y)$, that then $\nu(x^\sigma) = \nu(y^\sigma)$. An equivalent way to state the condition is $ \forall x \in K : \nu(x) = 0 \Rightarrow \nu(x^\sigma) = 0$ (to see this consider $\nu(xy^{-1})$ and $\nu(x^\sigma y^{-\sigma})$).

\subsection{Conclusion}\label{section:concl2}
In this section we will construct the root group labeling with epimorphism data $(u,\Theta,\theta,K,\Lambda,\nu,k)$, starting from a building $\Delta$ of type $(W,S)$ listed in Section~\ref{section:deffield} and $\phi: \Delta \to \Delta'$ an epimorphism of buildings.  

Choose a chamber $C$ of $\Delta$. Let $(u,\Theta, \theta)$ be the associated root group labeling in its standard form and $K$ its underlying field (see Section~\ref{section:high}). As we use the standard form there exist maps $\eta_s$ from $u(s)$ to the underlying field $K$ for every $s\in S$ by Section~\ref{section:rgl}. We also know by Section~\ref{section:subg} that the epimorphism gives rise to subgroups $w(s)\lhd v(s) \leq u(s)$. From the list of possible rank 2 residues it follows that there is at least one rank 2 residue containing $C$ to which we can apply the results of either Section~\ref{section:psl} or~\ref{section:ortho}. This yields a surjective valuation $\nu$ of the field $K$ to an ordered abelian group $\Lambda$ such that for a certain $t \in S$ and a $k(t) \in \Lambda$ one has that
\begin{align*}
v(t) &= \{g \in u(t) \vert  \nu(\eta_t(g)) \geq k(t) \}, \\
w(t) &= \{g \in u(t) \vert  \nu(\eta_t(g)) > k(t) \}.
\end{align*}
Repeated applications of Section~\ref{section:imp} (together with Remark~\ref{rem:comp}) to the root group sequences in the root group labeling $(u,\Theta, \theta)$ then imply that there exists constants  $k(s) \in \Lambda$ for each $s \in S$ such that 
\begin{align*}
v(s) &= \{g \in u(s) \vert  \nu(\eta_s(g)) \geq k(s) \}, \\
w(s) &= \{g \in u(s) \vert  \nu(\eta_s(g)) > k(s) \}.
\end{align*}
The compatibility conditions derived in Section~\ref{section:cc} yield conditions on these constants $k(s)$ and the algebraic structures defining $\Delta$. Note that the condition $k(s) \in \nu(\eta_s(u(s)))$ from Section~\ref{section:rgl} is indeed satisfied, this as the image of the epimorphism is not a weak building.

This concludes the description of the root group labeling with epimorphism data $(u,\Theta,\theta,K,\Lambda,\nu,k)$ arising from the epimorphism $\phi$ and hence the proof of the second Main Result. Note that this information determines the epimorphism by Corollary~\ref{cor:rigid2}.


\section{Proof of the third Main Result and the Main Corollaries}\label{section:construct}

The goal of this section is to construct an epimorphism given a root group labeling with epimorphism data $(u,\Theta,\theta,K,\Lambda,\nu,k)$ satisfying the compatibility conditions from Section~\ref{section:cc}. 

We will achieve this in multiple steps, generalizing the possible rank of the valuation $\nu$ at each step. We start with rank zero, then rank one, finite rank and eventually arbitrary rank, after which the proof of the third Main Result and the Main Corollaries will be concluded in Section~\ref{section:conclusion}. The second and third steps for generalized quadrangles of exceptional type will be handled separately.

Let $\Pi$ be the spherical building associated to the root group labeling $(u,\Theta,\theta)$. 

\subsection{Valuations of rank zero}\label{section:zero}
A rank zero valuation is trivial, so its value group consists only of the identity element. The corresponding epimorphism of $\Pi$ is the identity isomorphism from $\Pi$ to itself.

Note that one can view this isomorphism as an affine epimorphism (so arising from an $\R$-building), by considering the Euclidean cone over $\Pi$ as an $\R$-building (see for instance~\cite[Ex. 7.3]{Kra-Wei:*}).

\subsection{Valuations of rank one}\label{section:oner}
A rank one valuation is a valuation to the real numbers. Because of this the compatibility conditions will imply stronger restrictions. We start with the  octogonal systems in the following lemma.
\begin{lemma}
Let $K$ be a field equipped with a rank one valuation $\nu: K \to \R \cup \{\infty\}$. If for every $x \in K$ the implication $ \nu(x) = 0 \Rightarrow \nu(x^\sigma) = 0$ holds, where $\sigma$ is a Tits endomorphism of $K$, then $\nu(x^\sigma) = \sqrt{2} \nu(x)$.
\end{lemma}
\proof
Note that if $\nu(x) > 0$, then $\nu(x+1)$ and $\nu((x+1)^\sigma) = \nu(x+1)$ are both zero. So $\nu(x) \geq 0$ implies $\nu(x^\sigma) \geq 0$. This implies that $\nu$ and $\nu \circ \omega$ are equivalent valuations. By~\cite[Thm. 1.1.1]{Art:67} and the valuation $\nu$ being of rank one, this implies that there exists a positive real number $\gamma$ such that $\nu(x^\sigma) = \gamma\cdot  \nu(x)$ for all $x \in K^*$. The endomorphism $\sigma$ is a Tits endomorphism of a field of characteristic 2, so $(x^\sigma)^\sigma = x^2$. Applying $\nu$ to both sides of this equation yields that $\gamma$ equals $\sqrt{2}$. \qed

We now show what happens to the compatibility conditions consisting of inequalities, using quadratic forms as example. The compatibility condition here is the inequality
\[
\forall u,v \in L_0: \nu(q(u)),\nu(q(v)) \geq k \Rightarrow \nu(f(u,v)) \geq k.
\]
By replacing $u$ and $v$ by the scalar products $au$ and $a'v$ where $a,a' \in K$ we obtain
\begin{align*}
\forall u,v \in L_0: \nu(q(u)) \geq k - 2\nu(a), \nu(q(v)) \geq  k & -2 \nu(a')  \\
\Rightarrow  \nu(f(u,v)) & \geq k - \nu(a) -\nu(a').
\end{align*}

As $\nu$ is a valuation to the real numbers this implies an inequality of the following form where $C \in \R$ is a positive constant depending on the image of $\nu$ in $\R$.
\begin{align*}
\forall u,v \in L_0^* &:  \nu(f(u,v)) + C \geq (\nu(q(u)) + \nu(q(v)))/2 , 
\end{align*}
If the image of $\nu$ is dense in $\R$, then $C$ can be set 0. If it is $\mathbb{Z}$, then one can set $C$ to be 3. In~\cite[Prop. 19.4]{Wei:09} it is proven that this inequality with $C=0$ is equivalent to a condition for the completion of the quadratic form with respect to $\nu$. The part of the proof of this proposition that derives this condition from the inequality still holds if one only assumes the weaker inequality stated above. Hence the above inequality is equivalent to 
\begin{align*}
\forall u,v \in L_0^* &:  \nu(f(u,v))  \geq (\nu(q(u)) + \nu(q(v)))/2 .
\end{align*}
A similar reasoning for the other inequalities occurring in Moufang polygons which are not quadrangles of exceptional type is possible using~\cite[Prop. 24.9, 25.5 and 21.36]{Wei:09} (for the pseudo-quadratic form case, the substitution $(v,s)$ by $(va,a^\sigma s a)$ with $a \in K$ takes the role of the scalar product in the above discussion).

With the extra conditions we derived here one satisfies exactly the conditions (which can be found in~\cite{Hit-Kra-Wei:10} and~\cite{Wei:09}) for the existence of a valuation of the root data of the spherical building. By~\cite[Thm. 3]{Tit:86}, this implies the existence of an $\R$-building with the Moufang building $\Pi$ at infinity, corresponding to the valuation $\nu$ (see also~\cite[Thm. 14.47]{Wei:09} for the discrete case). Using the theory of $\R$-buildings one can obtain a canonical epimorphism of the spherical building at infinity (being $\Pi$) to a residue such that its root group labeling with epimorphism data is exactly the one we started with.

The fact that our compatibility conditions only concern the rank 2 residues is reflected in the result~\cite[Thm. 16.14]{Wei:09} on the existence of $\R$-buildings.

\begin{rem} \rm
The inequalities we derived in this section are not generally true if one leaves the rank one case. A consequence is that one cannot use the theory of $\Lambda$-buildings, which is the natural generalization of affine buildings for arbitrary valuations, to construct the epimorphism in one step. More information on $\Lambda$-buildings can be found in~\cite{Ben:90} and~\cite{Hit:09}.
\end{rem}


\subsection{Valuations of finite rank} \label{section:frank}
An abelian ordered group $\Lambda$ of rank $t$ can be embedded as a subgroup in the lexicographically ordered group $\oplus_{j=1}^t \R$ by Hahn's embedding theorem (see~\cite{Hah:07}). Using this presentation one can define an epimorphism $e: \Lambda \rightarrow \R : (a_1, \dots, a_t) \mapsto a_1$ of ordered abelian groups. Denote the kernel of this epimorphism by $\Lambda_0$. The function $\nu' := e \circ \nu$ is then a valuation of $K$ of rank one.

The claim is now that the compatibility conditions are satisfied for the valuation $\nu'$ as well. We again illustrate this with the octogonal sets and quadratic forms as examples.

For octogonal sets we have to prove that for $x\in K$ one has that $\nu(x) \in \Lambda_0 \Rightarrow \nu(x^\sigma) \in \Lambda_0$ given $\nu(x) = 0 \Rightarrow \nu(x^\sigma) = 0$. Suppose that this is not the case for a certain $x$, so $\nu(x) \in \Lambda_0$ and $\nu(x^\sigma) \notin \Lambda_0$. Note that $\nu(x^2) \in \Lambda_0$ as $\Lambda_0$ is a subgroup of $\Lambda$. Without loss of generality one may additionally assume that $\nu(x) < \nu(x^\sigma)$ (otherwise one can consider $x^{-1})$. In particular this implies that $\nu(1+x^2x^{-\sigma}) < 0$ (as $\Lambda_0$ is a convex subgroup).
 One also has $\nu(1+x^\sigma x^{-1}) = 0$, so $\nu((1+x^\sigma x^{-1})^\sigma)=\nu(1+ x^2 x^{-\sigma})  =0$, which is a contradiction. We conclude that $\nu(x) \in \Lambda_0 \Rightarrow \nu(x^\sigma) \in \Lambda_0$.

For quadratic forms, note that the inequality  
\[
\forall u,v \in L_0: \nu(q(u)),\nu(q(v)) \geq k \Rightarrow \nu(f(u,v)) \geq k.
\]
implies, by replacing $u$ and $v$ by scalar products $lu$ and $l'v$ with $l$ and $l'$ elements of the field with valuation in $\ker(e)$, that
\[
\forall u,v \in L_0: \nu(q(u)) + 2 \nu(l) ,\nu(q(v)) + 2 \nu(l') \geq k \Rightarrow \nu(f(u,v)) \geq k .
\]
As the valuation of $l$ and $l'$ can be chosen to be arbitrarily low in $\ker(e)$ when we apply the epimorphism $e$ at every place, we obtain that
\begin{equation}\label{eq:ineq}
\forall u,v \in L_0: \nu'(q(u)),\nu'(q(v)) \geq e(k) \Rightarrow \nu'(f(u,v)) \geq e(k) .
\end{equation}
We can now apply the results of the previous section to the rank one valuation $\nu'$, and obtain an epimorphism of the spherical building to some other spherical building which is defined over the residue field of $K_{\nu'}$. The valuation $\nu$ of $K$ allows us to define a rank $t-1$ valuation $\bar{\nu}$ of $K_{\nu'}$. On this new spherical building we can repeat the procedure until we have constructed the desired epimorphism, provided we can show the compatibility conditions for this new situation. We will be able to do this except for the Moufang quadrangles of exceptional type, for which there is no description (yet) of the possible residues. This is why we postpone this case to the next section.

For compatibility conditions involving an involution or Tits endomorphism $\sigma$ it is clear that the conditions stay true for a residue field. We will describe the behavior of conditions involving inequalities with the example of quadratic forms. As $\nu'$ is a rank one valuation, Equation~\ref{eq:ineq} implies by Section~\ref{section:oner} that
\[
\forall u,v \in L_0^* :  \nu'(f(u,v))  \geq (\nu'(q(u)) + \nu'(q(v)))/2 .
\]
The residue will be again a quadrangle of quadratic form type, where the quadratic space $\bar{L}$ is the quotient $\{v \in L_0 \vert \nu(q(v)) \geq e(k) \} / \{v \in L_0 \vert \nu(q(v)) > e(k) \} $ on which the residue field $K_{\nu'} $ acts naturally and for which the function $\bar{q}: \bar{L} \rightarrow K_{\nu'}: \bar{v} \mapsto \overline{q(v)/t}$ (where $\bar{.}$ indicated the natural map into $\bar{L}$ or $K_{\nu'}$, and $t \in K$ is such that $\nu'(t) =e(k)$) is an anisotropic quadratic function. See~\cite[Def. 19.33]{Wei:09} for more details to this construction. The original compatibility condition applied to $\bar{L}$ yields (keeping in mind the previous inequality to show independence of choice of representators)
\[
\forall u,v \in \bar{L}: \bar{\nu}(\bar{q}(\bar{u})),\nu(\bar{q}(\bar{v})) > k - \nu(t) \Rightarrow \nu(\bar{f}(\bar{u},\bar{v})) \geq k- \nu(t),
\]
where $\bar{f}$ is the bilinear form associated to $\bar{q}$. Hence we obtained a compatibility condition for the residue and we can continue with the construction of the epimorphism. For other types an analogous treatment is possible (see~\cite[Def. 24.50 and 25.28]{Wei:09} for detailed descriptions of the residues).

\subsection{Valuations of finite rank and quadrangles of exceptional type}\label{section:frank2}

In this section we handle Moufang quadrangles of exceptional type $\mathsf{F}_4$ or $\mathsf{E}_i$ ($i=6,7,8$). Combining the valuation $\nu$ of finite rank on the underlying field to the ordered abelian group $\Lambda$  with the norm functions $\eta_r$ listed in~\ref{section:imp}, we obtain maps $\phi_r: U_r \rightarrow \Lambda$ ($r \in \{1,2,3,4\}$). We are now interested in the interaction between these functions and the action of the $\mu$-maps of elements of $U_1$ and $U_4$. One observes using the relations in~\cite[32.10-11]{Tit-Wei:02} and~\cite[Prop. 21.10 and 22.4]{Wei:09} that (with $u_r \in U_r$)
\begin{align}
\phi_4(u_2^{\mu_1(u_1)}) &= \phi_2(u_2) - \phi_1(u_1), \notag \\
\phi_2(u_4^{\mu_1(u_1)}) &= \phi_4(u_4) + \phi_1(u_1),  \label{eq:mu1} \\
\phi_3(u_1^{\mu_4(u_4)}) &= \phi_1(u_1) + 2 \phi_4(u_4),  \notag \\
\phi_1(u_3^{\mu_4(u_4)}) &= \phi_3(u_3) - 2 \phi_4(u_4). \notag
\end{align}
Other identities are not straightforward to obtain. One can also derive that if $[u_1,u_4] = u_2 u_3$, then 
\begin{align*}
\phi_2(u_2) &= \phi_1(u_1) + \phi_4(u_4), \\
\phi_3(u_3) & = \phi_1(u_1) + 2 \phi_4(u_4). 
\end{align*}
One can use this in a reasoning similar to~\cite[Prop. 15.25]{Wei:09} (which makes use of the fact that double $\mu$-actions maps the root groups to themselves) obtaining (where $u_r, v_r, w_r \in U_r$)
\begin{align}
\phi_1(u_1^{\mu_1(v_1) \mu_1(w_1) }) &= \phi_1(u_1) - 2 \phi_1(v_1) + 2 \phi_1(w_1) \notag, \\
\phi_4(u_4^{\mu_4(v_4) \mu_4(w_4)}) &= \phi_4(u_4) - 2 \phi_4(v_4) + 2 \phi_4(w_4), \label{eq:mu} \\
\phi_3(u_3^{\mu_1(v_1) } ) &= \phi_3(u_3), \notag \\
\phi_2(u_2^{\mu_4(v_4)  }) &= \phi_2(u_2) . \notag
\end{align}
The compatibility conditions now are 
\begin{align}
\{u_1 \in U_1 \vert \phi(u_1) \geq k\} &\mbox{ is a subgroup}, \notag \\
\{u_1 \in U_1 \vert \phi(u_1) > k\} &\mbox{ is a subgroup}, \notag \\
\{u_4 \in U_4 \vert \phi(u_4) \geq l\} &\mbox{ is a subgroup}, \notag  \\
\{u_4 \in U_4 \vert \phi(u_4) > l\} &\mbox{ is a subgroup}, \notag \\
\phi_1(u_1) = k, \phi_3(u_3) >  k+2l, [u_1,u_3] &= u_2 \Rightarrow \phi_2(u_2) >k+l,\notag  \\
\phi_1(u_1) > k, \phi_3(u_3) =  k+2l, [u_1,u_3] &= u_2 \Rightarrow \phi_2(u_2) >k+l, \label{eq:except} \\
\phi_1(u_1) = k, \phi_3(u_3) =  k+2l, [u_1,u_3] &= u_2 \Rightarrow \phi_2(u_2) \geq k+l,\notag  \\
\phi_2(u_2) = k+l, \phi_4(u_4) > l, [u_2,u_4] &= u_3 \Rightarrow \phi_3(u_3) > k+2l,\notag  \\
\phi_2(u_2) > k+l, \phi_4(u_4) = l, [u_2,u_4] &= u_3 \Rightarrow \phi_3(u_3) > k+2l, \notag \\
\phi_2(u_2) = k+l, \phi_4(u_4) = l, [u_2,u_4] &= u_3 \Rightarrow \phi_3(u_3) \geq  k+2l.\notag 
\end{align}
From now on we keep in mind only the derived identities and inequalities and `forget' that we were dealing with an exceptional type case. The approach for constructing the epimorphism resembles the one from the previous two sections, but works in a more implicit way. As $\Lambda$ is of finite rank, one can find an ordered abelian group epimorphism $e : \Lambda \to \R$, which we compose with the maps $\phi_r$ to obtain maps $\phi'_r: U_r \to \R$ ($r \in \{1,2,3,4\}$). We claim that these maps satisfy the inequalities
\begin{align*}
[u_1,u_3] = u_2 \Rightarrow \phi'_2(u_2) \geq (\phi'_1(u_1) + \phi'_3(u_3))/2, \\
[u_2,u_4] = u_3 \Rightarrow \phi'_3(u_3) \geq (\phi'_2(u_2) + \phi'_4(u_4))/2, 
\end{align*}
and one has that 
\begin{align*}
\{u_1 \in U_1 \vert \phi'_1(u_1) \geq t \}, \\ 
\{u_4 \in U_4 \vert \phi'_4(u_4) \geq t \},
\end{align*}
with $t \in \R$ are subgroups. One can prove this by using Equations~\ref{eq:mu1} and~\ref{eq:mu} which allow us to add multiples of two to the constants $k$ and $l$ in the equalities of~\ref{eq:except}, . All of this implies that the maps $\phi'_r$ form a viable partial valuation in the sense of~\cite[Def. 15.4]{Wei:09}, and hence give rise to an $\R$-building (by adapting~\cite[Thm. 15.21]{Wei:09} to the non-discrete case) with the desired first epimorphism to a certain residue. The root groups $\overline{U}_r$ with $r\in \{1,2,3,4\}$ are given by the quotient $\{u_r \in U_r \vert \phi'_r(u_r) \geq 0 \} / \{u_r \in U_r \vert \phi'_r(u_r) > 0 \}$. On these root groups one can define in a natural way functions $\overline{\phi}_r: \overline{U}_r \rightarrow \ker e$, which inherit the same identities and inequalities as derived for the functions $\phi_r$. Hence we are back at our starting point and can apply recursion to obtain the desired epimorphism.
 


\subsection{Valuations of arbitrary rank}
In this final case we will construct an epimorphism of the spherical building $\Pi$ of type $(W,S)$ starting from the root group labeling with epimorphism data $(u,\Theta,\theta,K,\Lambda,\nu,k)$ not assuming any conditions on the rank of the valuation on the underlying field $K$ anymore. We will first prove this under the extra assumption that the (pseudo-)quadratic form appearing (if applicable) is finite-dimensional. The infinite-dimensional case will follow similarly after this is established. Let $\cC$ be the set of chambers of $\Pi$.

If a subfield $K_1$ of $K$ contains all the necessary constants in $K$ to ensure that the (pseudo-)quadratic form, hexagonal system and etc. are defined (note that the finite-dimensionality assumption assures that these constants are finite in number), and if $K_1$ is closed under the involution or Tits endomorphism $\sigma$ where applicable, then we call it a \emph{viable subfield}.

\begin{lemma}\label{lemma:fgviable}
Each finitely generated subfield of $K$ can be extended to a finitely generated viable subfield.
\end{lemma}
\proof
This holds as we only need to a finite number of generators: the constants to ensure that the forms are defined, and then images of the generators under the involution or Tits endomorphism $\sigma$ if applicable. \qed

Recall that the groups $u(s)$ from the root group labeling $(u,\Theta,\theta)$ of $\Pi$ correspond to the panels containing a certain chamber $C$ of $\Pi$, see Section~\ref{section:high}. If a subfield $K_1$ is viable, then we can restrict these groups to groups $u_{K_1}(s)$, which correspond to a subbuilding $\Pi(K_1)$ of $\Pi$ containing $C$ and with underlying field $K_1$. Let $\cC(K_1)$ be the set of chambers of $\Pi(K_1)$.

We now return to the problem of constructing an epimorphism. Let $\widetilde{K}$ be a finitely generated viable subfield of $K$ such that $k(s) \in \nu(\eta_s(u_{\widetilde{K}}(s)))$ for every $s \in S$ (compare this with the condition on the $k(s)$ in Section~\ref{section:rgl}). Such a field exists by Lemma~\ref{lemma:fgviable} as this only requires one group element for each group $u(s)$, so we can make sure that $\widetilde{K}$ satisfies this condition by assuring it contains a certain finite subset of elements of $K$.

Let $\Psi$ be the set of finitely generated viable subfields $K_1$ of $K$ containing $\widetilde{K}$. For each such $K_1 \in \Psi$ we have that  $k(s) \in \nu(\eta_s(u_{K_1}(s)))$ for every $s \in S$. Note that $\Psi$ is closed under intersection (this uses the fact that a subfield of a finitely generated field is again finitely generated).

An important observation is that for every $K_1 \in \Psi$ the natural restriction of the root group labeling with epimorphism data to this subfield again yields a root group labeling with epimorphism data. As a finitely generated field is of finite transcendence degree and Abhyankar's inequality (\cite[Lem. 1]{Abh:56})) states that the rank of a valuation on a field is at most its transcendence degree plus one, the restriction of the valuation $\nu$ to $K_1$ is of finite rank. So we can apply the results of Sections~\ref{section:frank} and~\ref{section:frank2} to construct an epimorphism $\varphi(K_1)$ from $\Pi(K_1)$ to some spherical building $\Pi'(K_1)$ with set of chambers $\cC'(K_1)$. The next lemma shows how the different epimorphisms are intertwined.

\begin{lemma}\label{lemma:psi_independent}
Let $K_1, K_2$ be two subfields in $\Psi$, suppose that $D$ and $E$ are two chambers contained in the intersection of the sets of chambers $\cC(K_1)$ and $\cC(K_2)$. Then the Weyl distance between the pairs of images $D^{\varphi(K_1)}$ and $E^{\varphi(K_1)}$ in $\Pi'(K_1)$ on one hand, and $D^{\varphi(K_2)}$ and $E^{\varphi(K_2)}$ in $\Pi'(K_2)$ on the other hand are the same.
\end{lemma}
\proof
As noted before, the subfield $K_3 := K_1 \cap K_2$ is also in $\Psi$. The previous paragraph guarantees the existence of three epimorphisms of $\Pi(K_3)$, namely the epimorphism $\varphi(K_3)$ and the restrictions of the epimorphisms $\varphi(K_1)$ and $\varphi(K_2)$ to $\Pi(K_3)$. However these only differ by an isomorphism of the image by Corollary~\ref{cor:rigid2} as the subgroups arising from these epimorphisms of $\Pi(K_3)$ are identical in each case. This proves the lemma.
\qed

The following lemma, which is the geometric counterpart of Lemma~\ref{lemma:fgviable} is crucial in our construction.
\begin{lemma}\label{lemma:sheaf}
Let $\Omega$ be a finite set of chambers of $\Pi$, then there is a field $K_1$ in $\Psi$ such that the subbuilding $\Pi(K_1)$ contains $\Omega$.
\end{lemma}
\proof
For a chamber in $\Omega$ to be contained in $\cC(K_1)$ it is sufficient that a finite amount of root elations of certain $U_\alpha$ are contained in $U_\alpha(K_1)$. The set of all the corresponding field elements in $K$ generates a field which can be extended to a finitely generated viable subfield by Lemma~\ref{lemma:fgviable}. \qed

We now use this lemma to define a function $\bar\delta$ between pairs of chambers of $\Pi$ which we  call the \emph{relative Weyl distance}. Let $(D,E)$ be such a pair, then Lemma~\ref{lemma:sheaf} implies that there exists a field $K_1 \in \psi$ such that $\Pi(K_1)$ contains both $D$ and $E$. We then set $\bar\delta(D,E)$ to be the Weyl distance between $D^{\varphi(K_1)}$ and $E^{\varphi(K_1)}$. By Lemma~\ref{lemma:psi_independent} this is independent of the choice of $K_1$. Remark that $\bar\delta(D,E) = \bar\delta(E,D)^{-1}$ and $\bar\delta(D,D)=1$.

\begin{lemma}\label{lemma:equivclass}
For any three chambers $D$, $E$ and $F$ of $\Pi$, if $\bar\delta(D,E) =1$, then $\bar\delta(D,F) = \bar\delta(E,F)$.  
\end{lemma}
\proof By Lemma~\ref{lemma:sheaf} there exists a subfield $K_1$ in $\Psi$ such that $\Pi(K_1)$ contains the three chambers. As $\bar\delta(D,E) =1$, one has that $D^{\varphi(K_1)}=E^{\varphi(K_1)}$, which on its turn implies that the Weyl distance in $\Pi'(K_1)$ between $D^{\varphi(K_1)}$ and $F^{\varphi(K_1)}$ is the same as between $E^{\varphi(K_1)}$ and $F^{\varphi(K_1)}$. By construction of $\bar\delta$  we conclude that $\bar\delta(D,F) = \bar\delta(E,F)$.  \qed

In particular it follows that being at trivial relative Weyl distance is an equivalence relation. Let $\bar\cC$ be the set of equivalence classes. Let $\varphi$ be the associated quotient map from $\cC$ to $\bar\cC$. By Lemma~\ref{lemma:equivclass} the relative Weyl distance between two chambers does not change if one replaces one chamber with an equivalent one, so one can extend the definition of the relative Weyl distance $\bar\delta$ to equivalence classes. We now claim that $\bar\Pi :=(\bar\cC,\bar\delta)$ is a building of type $(W,S)$ and that $\varphi$ is the desired epimorphism. We prove this in what follows by approximating $\bar\Pi$ by the buildings $\Pi'(K_1)$ with $K_1$ in $\Psi$.

\begin{lemma}\label{lemma:localiso}
For every subfield $K_1$ in $\Psi$, there exists a unique injective map $\bar\varphi(K_1)$ from $\cC'(K_1)$ to $\bar\cC$ such that the following diagram commutes (where the top arrow is the inclusion map).
$$\xymatrix{ \cC(K_1) \ar[d]_{\varphi(K_1)}  \ar[rr]  & & \cC \ar[d]^\varphi  \\
\cC'(K_1)  \ar[rr]_{\bar\varphi(K_1)} & &  \bar\cC }$$
Moreover for every two chambers $D,E \in \cC'(K_1)$ the Weyl distance between $D$ and $E$ equals $\bar\delta(D^{\bar\varphi(K_1)}, E^{\bar\varphi(K_1)})$. The restriction of $\bar\Pi$ to the image of $\bar\varphi(K_1)$ is a building of type $(W,S)$ isomorphic to $\Pi'(K_1)$ via the map $\bar\varphi(K_1)$.
\end{lemma}
\proof
There is a unique such injection as for every chamber $F \in \cC(K_1)$ the preimages $(F^{\varphi(K1)})^{\varphi(K1)^{-1}}$ and $(F^{\varphi})^{\varphi^{-1}}$ are identical. The statement involving the Weyl distances follows from the construction of $\bar\delta$. This then directly implies that $\bar\varphi(K_1)$ forms an embedding of $\Pi'(K_1)$ in $\bar\Pi$.  \qed

\begin{lemma}\label{lemma:sheaf2}
Let $\Omega$ be a finite subset of $\bar\cC$, then there is a field $K_1$ in $\Psi$ such that $\cC'(K_1)^{\bar\varphi(K_1)}$ contains $\Omega$.
\end{lemma}
\proof
This follows from the combination of Lemmas~\ref{lemma:sheaf},~\ref{lemma:localiso} and the fact that $\varphi$ is a quotient map and hence surjective. \qed

\begin{prop}
 $\bar\Pi$ is a spherical building of type $(W,S)$ and $\varphi$ is the desired epimorphism.
\end{prop}
\proof
To verify that $\bar\Pi$ is indeed a building we need to check conditions (WD1)-(WD3) from Section~\ref{section:buildings}. Note that each of these conditions involve at most three elements in $\bar\cC$. So given such a  subset of at most three elements in $\bar\cC$ we can pick a field $K_1$ in $\Psi$ by Lemma~\ref{lemma:sheaf2} such that $\cC'(K_1)^{\bar\varphi(K_1)}$ contains it. As Lemma~\ref{lemma:localiso} guarantees that $\bar\Pi$ restricted to $\cC'(K_1)^{\bar\varphi(K_1)}$ is a spherical building of type $(W,S)$, the conditions (WD1)-(WD3) hold for our choice of elements. In particular conditions (WD1)-(WD3) therefore hold in general and $\bar\Pi$ is a spherical building of type $(W,S)$.
The surjective map $\varphi$ preserves $s$-equivalency for each $s\in S$ by its construction and hence is an isomorphism. From the construction of the epimorphisms $\varphi(K_1)$ with $K_1 \in \Psi$, their relation with the epimorphism $\varphi$ derived in Lemma~\ref{lemma:localiso} and the construction of the root group labeling with epimorphism data in Section~\ref{section:concl2}, one concludes that $\varphi$ has the same root group labeling with epimorphism data as we started with. \qed

In order to handle the cases where infinite-dimensional (pseudo-)quadratic forms occur, one can approximate these by considering the restrictions to  finite-dimensional sub vector spaces, in an analogous way as we approximated the case of a valuation of arbitrary rank by restricting to certain subfields where the valuation is of finite rank.

\subsection{Conclusion}\label{section:conclusion}
In Section~\ref{section:imp} we derived, given an epimorphism $\phi$ of a Moufang building of type listed in Section~\ref{section:deffield}, a root group labeling with epimorphism data and  some compatibility conditions which are satisfied by it. We then continued in Section~\ref{section:construct} starting from a root group labeling with epimorphism data satisfying these conditions, constructing an epimorphism $\psi$ from it. Corollary~\ref{cor:rigid2} now implies that if we started with the one obtained from $\phi$, then $\phi$ and $\psi$ only differ by isomorphisms (as in the statement of the third Main Result).

This concludes the proof of the third Main Result.  The first Main Corollary follows from Sections~\ref{section:frank} and~\ref{section:frank2}, the second from Sections~\ref{section:oner} and~\ref{section:frank2} combined with Proposition~\ref{prop:factor}.








\end{document}